\documentclass[11pt]{article}
\usepackage{amsmath, graphicx}
\usepackage{graphicx,color}
\usepackage{amsthm}

\usepackage{amssymb}
\newtheorem{theorem}{Theorem}[section]

\newtheorem{lemma}[theorem]{Lemma}

\theoremstyle{definition}
\newtheorem{remark}[theorem]{Remark}

\def\ds{\displaystyle}
\def\forall{\hbox{for all}~}

\def\bfv{{\bf v}}

\def\bfw{{\bf w}}

\def\bfb{{\bf b}}

\def\bfn{{\bf n}}

\def\ve{\varepsilon}

\def\D{{\cal D}}

\def\R{{\mathbb R}}
\def\CC{{\mathbb C}}

\def\div{\hbox{div}\,}

\def\vp{\varphi}
\def\sym{\hbox{sym}\,}
\def\skew{\hbox{\rm skew}\,}

\def\vs{\vskip 2em}

\def\v{\vskip 1em}

\def\avint{-\!\!\!\!\!\!\int}
\def\begi{\begin{itemize}}
\def\endi{\end{itemize}}
\def\C{{\cal C}}
\def\ov{\overline}
\def\Tilde{\widetilde}

\def\bega{\begin{array}}
\def\enda{\end{array}}

\def\bel{\begin{equation}\label}
\def\eeq{\end{equation}}
\def\sqr#1#2{\vbox{\hrule height .#2pt
\hbox{\vrule width .#2pt height #1pt \kern #1pt
\vrule width .#2pt}\hrule height .#2pt }}
\def\square{\sqr74}
\def\endproof{\hphantom{MM}\hfill\llap{$\square$}\goodbreak}

\begin{document}
\title{\bf  A Model of Controlled Growth}
\vs
\author{Alberto Bressan$^{(*)}$ and Marta Lewicka$^{(**)}$\\
\,
\\
\small (*)~Department of Mathematics, Penn State University, \\
\small  University Park, PA 16802, USA.\\
\small  (**) Department of Mathematics, University of Pittsburgh,\\
\small 301 Thackeray Hall, Pittsburgh, PA 15260, USA. \\
\,
\\
E-mails: bressan@math.psu.edu, lewicka@pitt.edu }

\maketitle

\begin{abstract} 
We consider a free boundary problem for a system of PDEs, modeling the 
growth of a biological tissue.  A morphogen, controlling volume growth,  is produced 
by specific cells and then diffused and absorbed throughout the domain.   
The geometric shape of the growing tissue is determined by the instantaneous minimization of an elastic 
deformation energy, subject to a constraint on the volumetric growth.
For an initial domain with $\C^{2,\alpha}$ boundary,
our main result establishes the local existence and uniqueness
of a classical solution, up to a rigid motion. 
\end{abstract}

\section{Introduction}
\label{sec:0}
\setcounter{equation}{0}

Aim of this paper is to analyze a system of PDEs on a 
variable domain, describing the growth of a biological tissue.
Motivated by \cite{BG, BM, BGM},
we consider a living tissue containing some ``signaling cells",
which produce morphogen (i.e.,~a growth-inducing chemical).
This morphogen diffuses throughout the tissue and is partially 
absorbed.  A ``chemical gradient" is thus created:  the concentration 
of morphogen is not uniform, being larger in regions
closer to the signaling cells. In turn, this variable concentration determines 
a different volumetric growth in different parts of the living tissue. 
This can provide a mechanism for controlling the 
growth of the domain toward a desired shape.   

As customary, we describe biological  
growth in terms of a vector field $\bfv(\cdot)$, determining 
the motion of single cells within the tissue. 
Calling $u(\cdot)$ the concentration of morphogen, the constraint on volumetric growth is expressed by
\bel{div}\div \bfv ~=~ g(u)\,,\eeq
where $g:\R \to \R_+$ is a (possibly nonlinear) response function, satisfying $g(0)=0$. 
At any given time $t$, the vector field $\bfv$ is then determined
(up to a rigid motion) by the requirement that it minimizes
a deformation energy, subject to the constraint (\ref{div}).
The model is closed by the assumption that signaling cells
are passively transported within the tissue. 

Calling $\Omega(t)$ the region occupied by the tissue at time $t$,
and $w(t,\cdot)$   the concentration of signaling cells, we prove that the above model yields a 
well posed initial value problem. More precisely, our main theorems show that, 
if the initial domain $\Omega(0)=\Omega_0$ has 
$\C^{2,\alpha}$ boundary  and if the initial concentration $w(0,\cdot)$
lies in the H\"older space $\C^{0,\alpha}(\Omega_0)$ 
for some $0<\alpha<1$, then the system of evolution equations
determining the growing domain  has a classical solution, 
locally in time. Moreover, this solution is unique up to rigid motions, and 
preserves the regularity of the initial data.

A wide literature is currently available on free boundary problems modeling set growth,
see for example \cite{BCLM, BEL, CL, E, L, PS}. 
A major goal of these studies has been the mathematical description of tumor growth
\cite{BF3, CF3, CE, CF, F, FR}.   Compared with earlier works,
our model has various new features. 
On one hand, it contains a  transport equation for the density 
of morphogen-producing cells.   By varying the location and concentration of these cells, 
one can study how different shapes are produced.   
Another fundamental difference is that in our model 
the velocity field $\bfv$ is found as the minimizer of an
elastic deformation energy involving the $L^2$ norm of the
symmetric gradient of $\bfv$. On the other hand, 
in free boundary problems modeling flow in porous media 
one minimizes the $L^2$ norm of the velocity field $\bfv$ itself 
(with suitable constraints).  As a consequence, 
while the solutions in \cite{BF3, CF3, CE, CF, F, FR} 
are unique, the solutions  that we presently construct
are uniquely determined only up to rigid motions.

The remainder of this paper is organized as follows.
In Section~\ref{sec:1} we introduce the basic model and collect
the main notation.  Section~\ref{sec:2} contains 
some geometric lemmas on the representation of a family of sets with 
sufficiently smooth boundary.

The heart of the matter is worked out in Section~\ref{sec:3},
where we construct approximate solutions 
by a time discretization algorithm.  At each time step,
the density $u(\cdot)$ of morphogen satisfies a linear elliptic
equation accounting for production, diffusion, and adsorption.
Existence and regularity of solutions follow from 
standard theory \cite{GT}. In turn, the existence of a vector field $\bfv(\cdot)$ satisfying
the divergence constraint (\ref{div}) and minimizing a suitable
elastic deformation energy is proved relying on Korn's inequality.
A careful analysis shows that the system of equations
determining this constrained minimizer is elliptic in the
sense of Agmon, Douglis, and Nirenberg. Thanks to the
Schauder type estimates proved in \cite{ADN}, we thus obtain
the crucial a-priori bound on the norm 
$\|\bfv\|_{\C^{2,\alpha}}$.   Finally, the density $w(\cdot)$
of signaling cells is updated in terms of a linear transport equation
with $\C^{2,\alpha}$ coefficients, providing an estimate on how 
the norm $\|w\|_{\C^{0,\alpha}}$ grows in time.
Section~\ref{sec:4} contains some additional estimates,
showing that our approximate solutions depend continuously on the initial data.

In Section~\ref{sec:5} we state and prove our first main result,
on the existence of classical solutions, locally in time.   
The uniqueness of these solutions, up to rigid motions, is then proved in  
Section~\ref{sec:6}.  Two simple examples, where the growing domain 
$\Omega(t)$ can be explicitly computed, are discussed in Section~\ref{sec:75}.

The last two sections contain some supplementary material.
In Section~\ref{sec:7} we reformulate the problem using Lagrangian coordinates.  
Namely, we show that the growth of the living tissue can be
described by  an evolution equation for the coefficients of a Riemann metric
tensor on a fixed domain. Finally, an extension of our basic model is proposed in
Section~\ref{sec:8}, where we derive a set of equations describing 
the growth of a 2-dimensional surface embedded in $\R^3$,  regarded as a thin elastic shell.

\section{The basic model} 
\label{sec:1}
\setcounter{equation}{0}
Let $\Omega(t)\subset \R^d$ be the region occupied by a living tissue at time $t$,
in a space of dimension $d$. Cases  $d=2$ or $d=3$ are the most
relevant, however we formulate and prove our results in the general
case of arbitrary dimension. 

Assume that a morphogen 
is produced by cells located within the tissue.   
Denote by $w(t,x)$ the density of these cells at time $t$ and at a point $x\in \Omega(t)$.
Calling $u=u(t,x)$ the concentration of morphogen,   
we shall assume that $u$ satisfies a linear diffusion-adsorption equation 
with Neumann boundary conditions:
$$\left\{\bega{ll}
 u_t =\Delta u - u + w\qquad\qquad &x\in \Omega(t), \cr
 \langle\nabla u, \bfn\rangle = 0\qquad\qquad & x\in\partial
 \Omega(t). \enda\right. $$
Since the time scale of chemical diffusion is much shorter than the time scale 
of tissue growth, at any given time $t$ the solution of the 
above problem will be very close to an equilibrium, described by the elliptic equation
\bel{2}\left\{\bega{ll}
\Delta u - u + w=0\qquad\qquad &x\in \Omega(t),\cr
\langle\nabla u,\bfn\rangle =0\qquad\qquad & x\in\partial \Omega(t).\enda\right.
\eeq
We observe that, for every $w\in L
^2(\Omega(t))$, the solution $u$ of (\ref{2}) provides the unique
minimizer of a quadratic functional over the space
$W^{1,2}(\Omega(t))$. Namely,  it solves the problem
\begin{equation}\label{2v} \tag{M}
\hbox{minimize:}\quad J(u)~\doteq~\int_{\Omega(t)} \Big({|\nabla u|^2
\over 2} + {u^2\over 2} - wu\Big)~\mbox{d}x.
\end{equation} 
 
Next, we need an equation describing motion of cells within the tissue.   This is
determined by the expansion caused by volume growth.
Call $\bfv= \bfv(t,x)$ the velocity of the cell located at $x\in \Omega(t)$ at time $t$.
In our model, at each time $t$, the vector field $\bfv(t,\cdot)$ is determined 
as the solution to the constrained minimization problem
\begin{equation}\label{3}\tag{E}
\hbox{minimize:}\quad E(\bfv)~\doteq~{1\over 2} \int_{\Omega(t)} 
|\sym \nabla \bfv|^2~\mbox{d}x \qquad 
\hbox{subject to:}\quad \div \bfv~=~g(u)\,.
\end{equation}
Notice that $E(\bfv)$ can be regarded as the elastic energy of an
infinitesimal deformation (displacement). Throughout the paper, we
assume that the function $g:\R\to [0, \infty)$ satisfies 
\bel{gprop}g\in \C^3(\R),\qquad g(0)=0,\qquad g', ~ g'', ~ g''' \mbox{ are uniformly bounded.}
\eeq

Finally, we assume that the morphogen-producing cells are passively 
transported within the tissue. The transport equation below is
supplemented by assigning an initial distribution of hormone-producing
cells on the initial domain:
\begin{equation}\label{5}\tag{H}\left\{\bega{ll}
 w_t + \div(w \bfv)=0 & \qquad\qquad x\in \Omega(t), \cr
 w(0,x)=w_0(x) & \qquad\qquad  x\in \Omega(0)=\Omega_0.
\enda\right.
\end{equation}

Notice that, as soon as the velocity field $\bfv$ is known, we can recover 
$\Omega(t)$ as the set reached at time $t$ by trajectories starting in $\Omega_0$.
More precisely:
\begin{equation}\label{6}\tag{G}
\Omega(t)=\bigg\{ x(t)\,;\quad x(0)=x_0\in \Omega_0 ~\mbox{ and }~
 x'(s) =\bfv(s, x(s)) ~~\forall s\in [0,t]\bigg\}.
\end{equation}

Summarizing, we have:
\begi
\item[(i)] The linear  elliptic equation (\ref{2}), describing the concentration of morphogen
$u$ over the set $\Omega(t)$, at each time $t\geq 0$.
For a given source term $w(t,\cdot)$, 
its solution $u(t,\cdot)$ 
provides the unique minimizer in (\ref{2v}).

\item[(ii)] A constrained minimization problem (\ref{3}), determining the velocity field $\bfv(t,\cdot)$
at each given time $t$, up to a rigid motion: translation + rotation.

\item[(iii)] The linear transport equation (\ref{5}), determining how the concentration of morphogen-producing cells
evolves in time.

\item[(iv)] The formula (\ref{6}), describing the growth of the domain
  $\Omega(t)$.
\endi

The main goal of our analysis is to prove that, 
  given
an initial set $\Omega_0$ and an initial density
$w_0(x)$ for $x\in \Omega_0$, the equations (\ref{2v}-\ref{3}-\ref{5}-\ref{6})  determine a unique evolution (at least locally in time),
up to a rigid motion that does not affect the shape of the growing domain.

\subsection{Notation}
Throughout this paper, by $'$ or $\frac{{\rm d}}{{\rm d}t}$ we denote a derivative w.r.t.~time
$t$, while $\nabla$ is the gradient w.r.t.~the space variable $x
=(x_1,\ldots,x_d)$.

Given  a bounded, open, simply connected set $\Omega\subset\R^d$,
its boundary is denoted by $\Sigma=\partial \Omega$, 
and its Lebesgue measure by $|\Omega|$.
We write $\bfn$ for the outer unit normal vector to $\Omega$ at
boundary points, while $T_P(\partial\Omega)$ is the space 
of tangent vectors to the boundary $\partial\Omega$ 
at the point $P$. The average value of a function $f$ on $\Omega$ is denoted by
$$\avint_\Omega f~\mbox{d}x \doteq \frac{1}{|\Omega|}\int_\Omega f~\mbox{d}x.$$
For any integer $k\geq 0$ and $\alpha\in (0,1)$,
by $\C^{k,\alpha}(\Omega)$ we mean the space of bounded
continuous functions whose derivatives up to order $k$ are H\"older
continuous on $\Omega$, with the exponent $\alpha$.  This is a Banach
space with the norm:
$$\|u\|_{\C^{k,\alpha}(\Omega)}~\doteq~\sum_{|\nu|\leq k}~
\sup_{x\in\Omega} |\nabla^\nu u(x)| + \sum_{|\nu|=k}~\sup_{x,y\in\Omega
,\, x\not= y} ~{|\nabla^\nu u(x)- \nabla^\nu u(y)|\over |x-y|^\alpha}\,.$$
Since every H\"older continuous function $u$ as above
admits a unique extension to the closure $\ov\Omega$,
we observe that the spaces $\C^{k,\alpha}(\Omega)$ and 
$\C^{k,\alpha}(\ov\Omega)$ can be identified.

Given a $d\times d$ matrix $A= [A_{ij}]_{i,j=1\ldots d}$, we denote
by $A^T=[A_{ji}]$ its transpose, and we set:
$$ \sym A~\doteq~{A+A^T\over 2}\,,\qquad\skew A~\doteq~{A-A^T\over 2}\,,$$
$$ \langle A : B\rangle ~\doteq~ \mbox{trace}(A^TB) \,,
\qquad |A|^2~\doteq~ \langle A:A\rangle = \sum_{i,j=1}^dA_{ij}^2.$$
The space of $d\times d$ skew-symmetric matrices is $so(d)$, and $I$
is the $d\times d$ identity matrix.  
%(we always think $d=2$ or $d=3$).

%Here the key issue is to provide estimates on how the solution to the constrained
%minimization problem (\ref{3})-(\ref{4}) depends on the data.   If this dependence is
%sufficiently regular, one can hope to construct a convergent sequence of approximations.

%\begin{figure}[htbp]
%\centering
 % \includegraphics[scale=0.4]{FIG/pg13.eps}
  %  \caption{{\small An example of set growth. 
  %The shaded areas denote the support    of $w$.}}
%\label{f:pg1}
%\end{figure}

\section{Some geometric lemmas}\label{sec:2}
\setcounter{equation}{0}

We say that $\Omega$ satisfies the uniform inner and outer sphere
condition when there exists $\rho>0$ such that, for every 
boundary point 
$x\in \Sigma$, we can find closed balls $B^{in}$ and $B^{out}$ 
of radii $R_{{in}}(x), R_{{out}}(x) \geq 2\rho$ satisfying $B^{in} \subset \ov \Omega$,
$B^{in}\cap\Sigma =\{x\}$ and $B^{out}\cap\ov\Omega =\{x\}$. Define the signed distance function:
$$\delta(x)~\doteq~\left\{ \bega{cl} 
~{\rm dist}(x,\Sigma)\qquad & ~~x\notin \Omega\cr
-{\rm dist}(x,\Sigma)\qquad & ~~x\in \Omega.\enda\right.$$
If $\Omega$ is smooth (i.e., it has a smooth boundary), 
then $\delta(\cdot)$ is also
smooth, when restricted to the open set
$$V_\rho~\doteq ~\{ x;~~{\rm dist} (x,\Sigma)<\rho\}.$$
Moreover,  for every $x\in V_\rho$ there exists a unique point
$\pi(x)\in\Sigma$  with $|\pi(x)-x|= {\rm dist}(x, \Sigma).$

\begin{figure}[htbp]
\centering
\includegraphics[scale=0.45]{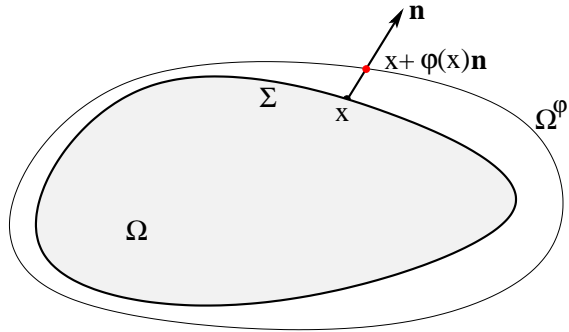}
    \caption{{\small The set $\Omega^\vp$ in (\ref{omp}),
    described in terms of the function $\varphi:\Sigma\to\R$.}}
\label{f:sg32}
\end{figure}

Every continuous map $\varphi:\Sigma\to (-\rho,\rho)$
determines then a bounded open set (see Fig. \ref{f:sg32}): 
\bel{omp}
\Omega^\varphi=\bigl\{ x\in\R^d;~~\delta(x)< \varphi(\pi(x))\bigr\} \quad
\mbox{with} \quad \partial\Omega^\varphi=\bigl\{y+\varphi(y)\,\bfn(y);~~y\in \Sigma\bigr\}.
\eeq
To measure the H\"older regularity of $\varphi$, we extend it
to  $V_\rho$ by  $\varphi(x)\doteq\varphi(\pi(x))$, and set:
\begin{equation}\label{definorm}
\|\varphi\|_{\mathcal{C}^{k, \alpha}} ~\doteq ~ \|\varphi\|_{\mathcal{C}^{k, \alpha}(V_\rho)}.
\end{equation}

By definition, $\Sigma\in \C^{k,\alpha}$ if the following holds.
For every $x\in \Sigma$ there exists an open  ball 
$B(x,r)$ and a homeomorphism $h:B(x,r)\to B(0,1)\subset \R^d$ such
that \vspace{-2mm}:
\begin{itemize}
\item[(i)] The map $h$ as well as its inverse $h^{-1}$ are $\C^{k,\alpha}$ regular.
\item[(ii)] $h\big(B(x,r)\cap\Omega\big) = B(0,1)\cap \bigl\{x\in\R^d;~x_1>0\bigr\}.$
\end{itemize}

\begin{lemma}\label{lem2.3}
Let $\Omega\subset \R^d$ be an open, bounded, simply connected and
smooth set, satisfying the uniform inner and outer sphere condition
with radius $2\rho >0$. Then, 
for every $\kappa>0$ there exists a constant $M$ such that the following holds.
If $\varphi:\Sigma\to (-\frac{\rho}{2}, \frac{\rho}{2})$ satisfies
$\|\varphi\|_{\C^{2,\alpha}}\leq \kappa$, then there exists a
homeomorphism $\Lambda: \Omega \to\Omega^\varphi$ satisfying the bounds:
\bel{Cka}
\|\Lambda\|_{\C^{2,\alpha} (\Omega)}\leq M,\qquad\quad 
\|\Lambda^{-1}\|_{\C^{2,\alpha}(\Omega^{\varphi})} \leq M.
\eeq
\end{lemma}

{\bf Proof.} {\bf 1.}
Let $\sigma:\R\to\R$ be a $\C^\infty$ function such that
$\sigma(s)=0$ for $s\leq -{\rho }$, and $\sigma(x) = 1$ for $s\geq 0$, and moreover:
\bel{sprop}
0\leq \sigma'(s)\leq \frac{3}{2\rho} \quad ~\hbox{for all} \quad s\in\R.
\eeq
The homeomorphism $\Lambda:\Omega\to\Omega^\varphi$ is defined by setting:
$$\Lambda(x)=\left\{ \bega{cl} x\qquad &\hbox{if}\quad ~~\delta(x)\leq -\rho\\
x + \sigma(\delta(x))\varphi(x) \bfn(\pi(x))\qquad &\hbox{if}\quad -\rho<\delta(x)<0.
\enda\right.$$
It is easily seen that $\Lambda$ maps $\Omega$ onto $\Omega^\varphi$.
Since  $\Lambda$ coincides with identity 
on the set where $\delta(x)\leq -\rho$,
to estimate the $\C^{2,\alpha}$ norm of $\Lambda$ it suffices to
study what happens when $-\rho<\delta(x)<0$.   
On this latter set, the functions $\delta(x)$, $\sigma(\delta(x))$, $\bfn(\pi(x))$ 
have uniformly bounded derivatives up to any order.  By the definition  of $\Lambda$
we thus get the estimate:
$$\|\Lambda\|_{\C^{2,\alpha}(\Omega)}\leq C \big(1+ \|\varphi\|_{\C^{2,\alpha}}\big),$$
for a suitable constant $C$ depending only on $\Sigma$.

{\bf 2.}  In order to obtain a similar estimate for $\Lambda^{-1}$, it is enough to check
that $\det\nabla\Lambda$ has uniformly bounded inverse on
$\Omega$. Indeed, in this case, the $\C^{2, \alpha}$ norm of
$\Lambda^{-1}$ will be bounded by a polynomial in $\|\Lambda\|_{\C^{2,
  \alpha}(\Omega)}$ whose order and coefficients depend only on
$\Omega$ and $d$.

On the set where $\delta(x)\leq - \rho$, we have $\det\nabla \Lambda =
1$. Let now $-\rho<\delta(x)<0$, and let $y=\pi(x)\in\Sigma$.
Let $U\subset \Sigma$ be a relatively open neighborhood of $y$, with coordinates
$(x_2,\ldots, x_d)$.  Then the map $x\mapsto (\delta(x), x_2, \ldots, x_d)$
provides a chart of the inverse image $\pi^{-1}(U)$.  
In these coordinates, $\Lambda$ has the form: 
$$\Tilde \Lambda (x_1,\ldots , x_d)~ =~ \big( x_1 + \sigma(x_1)\varphi(x),
x_2,\ldots , x_d\big).$$
In view of (\ref{sprop}) and the fact that $\varphi$ is independent of $x_1$, we thus conclude:
$$\det\nabla\Tilde\Lambda(x)~ = ~1+\sigma'(x_1) \phi(x)~\geq~
1-\frac{3}{2\rho}\frac{\rho}{2} ~=~ \frac{1}{4}.$$
The estimate (\ref{Cka}) now follows by covering the compact surface 
$\Sigma$ with finitely many coordinate charts and by noting that,
on each chart, $\det\nabla\Lambda$ is uniformly comparable with 
$\det\nabla\Tilde\Lambda$. 
\endproof

\begin{lemma} \label{lem2.4}
Let $\Omega_0\subset \R^d$ be an open, bounded and simply connected
set with $\C^{2,\alpha}$ boundary $\Sigma_0$, satisfying the uniform 
inner and outer sphere condition with radius $3\rho>0$.  
Then, for any $\ve_0>0$, there exists an open, bounded and simply
connected set $\Omega$ with $\C^\infty$ boundary $\Sigma$, 
satisfying the uniform inner and outer sphere condition with radius $2\rho$, 
and such that $\Omega_0 = \Omega^\varphi$ as in (\ref{omp}) for
some function $\varphi\in \C^{2,\alpha}(\Sigma)$ with:
\begin{equation}\label{cos}
|\varphi (x)|<\varepsilon_0 \quad \mbox{for all} \quad
x\in\Sigma. 
\end{equation}
\end{lemma}

{\bf Proof.} {\bf 1.} Let $\delta_0$ be the signed distance function
from $\Sigma_0$. By assumption, $\delta_0$ is $\C^2$ on the open neighborhood
$V_{0,3\rho}$ of $\Sigma_0$ with radius $3\rho$. 
We now consider the mollification
$\delta_\ve = \delta_0 * J_\ve$ with a standard mollifier $J_\ve$ in
$\R^d$. It is not restrictive to assume that $\ve\ll\ve_0\ll\rho$ 
and that
\begin{equation}\label{pom}
\|\delta_\ve - \delta_0\|_{\C^{2,\alpha}(V_{0,3\rho-\ve_0})} \leq C\ve.
\end{equation}
We claim that the set
$$\Omega  = \Omega_\ve~\doteq~\{ x\in\R^d;~~\delta_\ve(x)<0\}$$
satisfies the conclusions of the lemma, provided that  $\varepsilon>0$ is
chosen
sufficiently small. Since
$|\nabla \delta_0|=1$ in $V_{0, 3\delta}$, we note that:
$$|\nabla \delta_\ve(x)|\geq 1-{\ve_0\over 2}\quad \forall x\in V_{0, 3\rho}, \qquad
|\delta_\ve(x)|\leq {\ve_0\over 2}\quad \forall x\in \Sigma_0. $$
Now fix $x\in \Sigma_0$.   By the above estimates and
since $\delta_0\in\mathcal{C}^2$, we can find $y\in V_{0, \rho}$ such that
$$\delta_\ve(y)=0\qquad \hbox{and}\qquad |y-x|\leq {\ve_0\over
  2} \left( 1-{\ve_0\over 2}\right)^{-1}<\ve_0.$$ 
Consequently, every point $x\in \Sigma_0$ is at a distance less than $\ve_0$ from some
$y\in \Sigma_\ve =\partial\Omega_\ve$.
We conclude that the smooth set $\Omega= \Omega_{\ve}$ 
indeed satisfies $\Omega^\varphi = \Omega_0$ and the uniquely
determined function $\varphi$, given as the signed distance from
$\Sigma$, obeys (\ref{cos}) and it is $\C^{2,\alpha}$ regular.

\begin{figure}[htbp]
\centering
\includegraphics[scale=0.5]{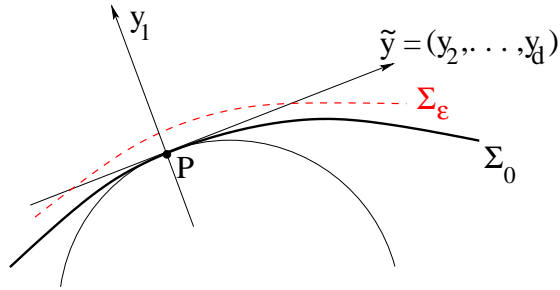}
    \caption{\small  Estimating the radius of curvature of the boundary 
    $\Sigma_\ve= \partial \Omega_\ve$ }
\label{f:sg51}
\end{figure}

{\bf 2.} We now check that $\Omega=\Omega_\ve$ satisfies the uniform inner
and outer sphere condition with radius $2\rho$. 
Fix any point $P\in \Sigma_0$.    On a neighborhood of $P$ 
we introduce an orthonormal frame of coordinates $(y_1,\ldots, y_d) =
(y_1,\tilde y)$ as in Fig. \ref{f:sg51}, where the $y_1$-axis is orthogonal to the surface 
$\Sigma_0$ at $P$.   In these local coordinates, the surfaces $\Sigma_0$, $\Sigma_\ve$  have the representations:
$$ \Sigma_0 = \{(y_1,\tilde y);~~y_1=\psi_0(\tilde y)\bigr\}, \qquad
\Sigma_\ve = \{(y_1,\tilde y);~~y_1=\psi_\ve(\tilde y)\bigr\},$$
with the variable $\tilde y$ ranging in some neighborhood of the origin  $U\subset \R^{d-1}$.

By construction we have $\frac{\partial\delta_0}{\partial y_1}(P) = 1$. 
Hence,  by possibly shrinking the neighborhood $U$, we can assume 
$\frac{\partial\delta_0}{\partial y_1}(\tilde y) \geq
\frac{1}{2}$ for every $\tilde y\in U$. By (\ref{pom}) we thus have
$\|\psi_\ve - \psi_0\|_{\C^0(U)}\leq C\ve$ and the implicit function theorem further implies the convergence
\bel{conv}
\|\psi_\ve-\psi_0\|_{\C^2(U)}~\to~0\qquad\qquad\hbox{as}\quad \ve\to 0\,.
\eeq

We now recall that the maximal curvature $\chi(\tilde y)$ of the graph of a
function $\psi: \R^{d-1}\to\R$ at a point $\tilde y$, equals the maximum of
the absolute values of the principal curvatures, i.e. the maximum of
the absolute values of the eigenvalues of the second fundamental form
$\Pi = (\nabla \psi)^T\nabla \mathbf{n}$. Since the second fundamental
forms of $\Sigma_0$ and $\Sigma_\ve$ satisfy: $\|\Pi_\ve
-\Pi_0\|_{\C^{0}(U)}\to 0$ as $\ve\to 0$ in virtue of (\ref{conv}),
and since for every $\tilde y\in U$ the assumption of the lemma gives:
$\chi_0(\tilde y)\leq\frac{1}{3\rho}$, it indeed follows that
$\chi_\ve(\tilde y) \leq \frac{1}{2\rho}$ for small $\ve>0$.

%More precisely, $\chi(\tilde y)$ may be obtained by taking the
%supremum over all directional curvatures, i.e. curvatures corresponding
%to all unit vectors $v\in\R^{d-1}$, at the point $(\psi(\tilde y), \tilde y)$:
%\bel{curvrad}
%\chi(\tilde y) = \sup_{|v|=1}{|\langle \nabla^2\psi (\tilde y) v, v\rangle|\over 
%\bigl(1+ \langle\nabla\psi (\tilde y), v\rangle^2\bigr)^{3/2}}.\eeq
%By assumption, the maximal curvature for $\psi_0$ satisfies:
%$\chi_0(\tilde y)\leq{1\over 3\rho}$  for all $\tilde y=(y_2,\ldots, y_d)$ in the domain of the chart.
%Thus, by (\ref{conv}) the maximal curvature for $\psi_\ve$ for all
%$\ve>0$ small enough, must also satisfy: $\chi_\ve(\tilde y) \leq{1\over 2\rho}$.
In turn, this yields an a-priori bound on the inner and outer curvature radii:
$$\min\bigl\{ R_{in}(\psi_\ve(\tilde y),\tilde y), ~ R_{out}(\psi_\ve(\tilde y), \tilde y)\bigr\}
~=~ {1\over \chi_\ve (\tilde y)}~ \geq~ 2\rho.$$
By covering the compact surface $\Sigma_0$ with neighborhoods of
finitely many points $P_1,\ldots, P_\nu$, and choosing
$\ve=\min\{\ve_1,\ldots,\ve_\nu\}$, the proof is achieved. \endproof

\section{Regularity estimates}
\label{sec:3}
\setcounter{equation}{0}

Given the initial data $w_0$ in (\ref{5}),  a local solution to the
system of equations (\ref{2v}-\ref{3}-\ref{5}-\ref{6}) 
will be constructed as a limit of approximations, obtained by 
discretizing time.

Fix a time step $\epsilon>0$ and let $t_k= k\epsilon$.  Assume that at time $t_k$ we are given the set
$\Omega_k=\Omega(t_k)$ and the scalar nonnegative function $w_k
=w(t_k,\cdot)$ on $\Omega_k$. Successive $\Omega_{k+1}= \Omega(t_{k+1})$ and $w_{k+1} =
w(t_{k+1}, \cdot)$ on $\Omega_{k+1}$ are obtained by the application of
the four steps below.

\begin{description}
\item[Step 1.]	Determine the density $u_k:\Omega_k\to\R$
by minimizing  (\ref{2v}) with $w=w_k$.  This implies that $u_k$ is the solution to
the elliptic problem (\ref{2}).
\item[Step 2.]	Determine the velocity field $\bfv_k:\Omega_k\to\R^d$ by solving 
the minimization problem (\ref{3}) on $\Omega_k$ subject to the
current constraint $\div \bfv_k= g(u_k)$.
The minimum is defined up to a rigid motion and we can single out 
a unique $\bfv_k$ by requiring that
\bel{vkno}\avint_{\Omega_k} \bfv_k~\mbox{d}x = 0\,,\qquad\qquad  \skew
\displaystyle{\avint_{\Omega_k} \nabla 
\bfv_k~\mbox{d}x =0.}\eeq 
\item[Step 3.]	Define the domain $\Omega_{k+1}$ by an approximation of (\ref{6}):
\bel{Ok1}
\Omega_{k+1}~\doteq~ \bigl\{x + \epsilon \bfv_k(x);~ x \in
\Omega_k\bigr\}.\eeq
\item[Step 4.] On the set $\Omega_{k+1}$, define the density $w_{k+1}$ implicitly by setting
\bel{wk1}
w_{k+1}( x + \epsilon \bfv_k(x))~ \doteq ~{w_k(x)\over  \det(I+\epsilon \nabla \bfv_k(x))}\,.
\eeq
Notice that (\ref{wk1}) is motivated by mass conservation:
 $w_{k+1}$ is the push-forward of the density $w_k$ through the map
$x\mapsto x+\epsilon \bfv_k(x)$.  The motivation for (\ref{wk1}) in
the continuous framework is given in Lemma \ref{motivation}.
\end{description}

Throughout the following, we assume that the initial domain
$\Omega_0\subset\R^d$ is open, bounded and simply connected,
with boundary $\Sigma_0\in \C^{2,\alpha}$, whereas the initial density
satisfies $w_0\in \C^{0,\alpha}(\Omega_0)$,  for some $0<\alpha<1$.
Moreover, the function $g\in \C^3(\R)$ satisfies (\ref{gprop}) unless
stated otherwise.

\subsection{Step 1: The elliptic equation for $u$}

\begin{lemma}\label{lemS1}
Let $\Omega\subset \R^d$  be an open, bounded and
simply connected set with $\C^{2,\alpha}$ boundary. Let $w \in
\mathcal{C}^{0,\alpha}(\Omega)$ be a nonnegative function. 
Then (\ref{2}) has a unique solution $u \in \mathcal{C}^{2,\alpha}(\Omega)$, which is
nonnegative and satisfies:
\bel{S1}
\|u\|_{\mathcal{C}^{2,\alpha}(\Omega)}~\leq ~C \|w\|_{\mathcal{C}^{0,\alpha}(\Omega)}.
\eeq
Further, for every constant $M>0$ and every domain $\Tilde\Omega$ for
which there exists a homeomorphism $\Lambda:\Omega\to\Tilde\Omega$ with
$\|\Lambda\|_{\C^{2,\alpha}(\Omega)},
\|\Lambda^{-1}\|_{\C^{2,\alpha}(\Tilde\Omega)}\leq M$, the 
corresponding bound  (\ref{S1})
is valid with a uniform constant $C$ that depends only on $M$ (in
addition to $\Omega$ and $\alpha$ that are given in the problem).
\end{lemma}
\v
{\bf Proof.} {\bf 1.} Existence and uniqueness of solutions to
(\ref{2}) follow from Theorem 6.31 in \cite{GT} 
(see also the remark at the
end of Chapter 6.7 in \cite{GT}). 
We now show the non-negativity of $u$. 
If $u$ is constant then $u=w\geq 0$. For
non-constant $u$, we invoke the maximum principle (Theorem 3.5
\cite{GT}) and conclude that the non-positive minimum of $u$ on
$\ov\Omega$ cannot be achieved in the interior $\Omega$. On the other hand, if
such minimum is achieved at some $x\in\partial\Omega$, then by Hopf's
lemma (see Lemma 3.4 in \cite{GT}), 
one must have $\langle \nabla u(x),
\mathbf{n}\rangle <0$, contradicting the boundary condition in 
(\ref{2}).
\v
{\bf 2.} Let now $\Lambda$ and $M$ be as in the statement of the
lemma. Let $\tilde u$ be the solution to (\ref{2}) on $\Tilde \Omega$,
for some $\tilde w \in \mathcal{C}^{0,\alpha}(%\ov
{\Tilde\Omega})$. Then the composition
$u=\tilde u\circ \Lambda\in \mathcal{C}^{2,\alpha}(%\ov
\Omega)$ provides the unique solution to the following boundary value problem:
\bel{pullback}\left\{\bega{ll}
\langle\nabla^2u : A\rangle + \langle \nabla u, \,\Delta
(\Lambda^{-1})\circ\Lambda\rangle - u ~=~ - \tilde w\circ \Lambda \qquad\qquad &x\in \Omega,\\[4mm]
\langle\nabla u, A \bfn\rangle~ =~0\qquad\qquad & x\in\partial \Omega.\enda\right.
\eeq
Here the matrix of coefficients $A$ is defined as
$$A(x)~ = ~\Big((\nabla \Lambda^{-1})(\nabla \Lambda^{-1})^T\Big)
(\Lambda(x)) = \Big((\nabla \Lambda)^T(\nabla \Lambda)\Big)^{-1}(x).$$
To derive the boundary condition, we used the
following formula which is valid for every invertible matrix: 
$(B\xi_1)\times (B\xi_2) = (\det B) B^{-1, T}(\xi_1\times \xi_2)$.  By
Theorem 6.30 in \cite{GT} we obtain the bound:
\begin{equation}\label{esti}
\|u\|_{\mathcal{C}^{2,\alpha}(\Omega)} ~\leq~ C
\left(\|u\|_{\mathcal{C}^{0,\alpha}(\Omega)} +  \|\tilde w\circ\Lambda\|_{\mathcal{C}^{0,\alpha}(\Omega)}\right),
\end{equation}
where the constant $C$ depends only on $\Omega$, $\alpha$ and on an upper
bound to the following quantities: $\|A\|_{\C^{1,\alpha}(\Omega)}$,
$\|\Delta(\Lambda^{-1})\circ\Lambda\|_{\mathcal{C}^{0,\alpha}(\Omega)}$
and the joint ellipticity and non-characteristic boundary constant $\kappa_\Lambda$. The defining
requirement for $\kappa_\Lambda$ is that:
$$\frac{1}{\kappa_\Lambda} |\xi|^2~\leq~ \langle A(x) \xi, \xi\rangle
~\leq ~{\kappa_\Lambda} |\xi|^2 \qquad \mbox{for all}~x\in %\ov
\Omega.$$ 
Hence we can simply take $\kappa_\Lambda = \|(\nabla \Lambda)^{-1}\|_{\C^0} + \|\nabla \Lambda\|_{\C^0} ^{2}$, 
confirming that the constant $C$ in (\ref{esti}) depends only on $M$.
\v
{\bf 3.} We now show that (\ref{esti}) can be improved to
\begin{equation}\label{esti2}
\|u\|_{\mathcal{C}^{2,\alpha}(\Omega)}~ \leq ~
C \|\tilde w\circ\Lambda\|_{\mathcal{C}^{0,\alpha}(\Omega)},
\end{equation}
for a possibly larger constant $C$, which still depends only on the bounding
constant $M$.  We argue by contradiction; assume there are sequences
of diffeomorphisms $\Lambda_n$ such that $\|\Lambda_n\|_{\C^{2,\alpha}},
\|\Lambda_n^{-1}\|_{\C^{2,\alpha}}\leq M$, and of solutions
$u_n\in\C^{2,\alpha}(%\ov
\Omega)$ to the problem (\ref{pullback}) with
some $\tilde w_n \in \C^{0,\alpha}(%\overline
{\Lambda_n(\Omega)})$, so that:
$$\|u_n\|_{\C^{2,\alpha}(\Omega)}=1 \qquad \mbox{and} \qquad \|\tilde
w_n\circ\Lambda_n\|_{\C^{0,\alpha}(\Omega)}\leq \frac{1}{n}.$$
Fix $\beta\in (0, \alpha)$. Passing to a subsequence if necessary, we
may assume that $\Lambda_n$ converge as $n\to\infty$ (together with their inverses)
in $\C^{2,\beta}(%\ov
\Omega)$ to some $\Lambda$, and that, likewise, $u_n$ converge to
$u$. The limit $u$ must then solve the problem (\ref{pullback}) with $\tilde
w = 0$. Thus $u=0$ and  $\|u_n\|_{\C^{0,\alpha}}$
converging to $0$ implies, in view of (\ref{esti}), that
$\|u_n\|_{\C^{2,\alpha}}$ converges to $0$ as well. This is a  contradiction
that achieves (\ref{esti2}). 

Noting that $\|\tilde u\|_{\C^{2,\alpha}}\leq C \|u\|_{\C^{2,\alpha}}$ and 
$\|\tilde w\circ\Lambda\|_{\C^{2,\alpha}}\leq C \|\tilde w\|_{\C^{2,\alpha}}$
with $C$ depending only on $M$, we see that (\ref{esti2}) yields
(\ref{S1}) on $\Tilde\Omega$.
\endproof

\subsection{Step 2: The elastic minimization problem for $\bfv$}

\begin{lemma}\label{lem1}
Let $\Omega\subset \R^d$  be an open, bounded and simply connected set with $\C^{2,\alpha}$ boundary.
Assume that  $u \in  W^{ 1,2}(\Omega, \R)$ and that
$g\in\mathcal{C}^1$ satisfy $g(0)=0$ with $g'$ bounded. Then the following holds.
\begi
\item[(i)] The minimization problem (\ref{3}) has a solution, which is unique up to rigid motions.
\item[(ii)] A vector field $\bfv\in W^{1,2}(\Omega,\R^d)$ is a
  minimizer of (\ref{3}) if and only if there exists $p\in L^2(\Omega, \R)$ such that
$(\bfv, p)$ solves:
\bel{EL}
\left\{ \bega{rll} \mathrm{div}(\mathrm{sym}\nabla \bfv - p I) & = ~0 \qquad  &x\in  \Omega,\\[2mm]
{\mathrm{div}}~\bfv& = ~g(u) &x\in  \Omega,\\[2mm]
({\mathrm{sym}} \nabla \bfv - p I) \bfn & = ~0 \qquad &x\in\partial\Omega.\enda\right.\eeq

\item[(iii)]  There exists a constant $C$, independent of $u$, such that
any $(\bfv,p)$ as above satisfies:
\bel{basic} 
\left\|\nabla \bfv - \mathrm{skew} \avint_\Omega \nabla \bfv~\mathrm{d}x
\right\|_{L^2(\Omega)} + 
\left\|p - \avint_\Omega p~\mathrm{d}x\right\|_{L^2(\Omega)}~\leq~ C \|u\|_{L^2(\Omega)}.\eeq
\endi
\end{lemma}

{\bf Proof.} {\bf 1.} Note that $g(u)\in W^{1,2}(\Omega, \R)$.
Existence in (i) follows by the direct method of Calculus of
Variations. Consider  a minimizing sequence $\bfv_n$.
By  Korn's and Poincar\'e's inequalities, we can replace  each $\bfv_n$ by a vector field of the form:
$$\tilde{\bfv}_n(x)~ =~ \bfv_n(x)
-(A_nx+\bfb_n),$$ where $A_n\in so(d)$ and  $\bfb_n\in {\R}^d$, so that
$\tilde\bfv_n\rightharpoonup \bfv$ weakly in $W^{1,2}$, up to a subsequence.  By the convexity of
the functional $E$, it is clear that the limit $\bfv$ is a minimizer.

To prove uniqueness, let $\bfv_1$ and $\bfv_2$ be two minimizers. Test the
minimization in (\ref{3}) in both $\bfv_1$ and $\bfv_2$ by the admissible divergence-free
perturbation field $\bfv_1-\bfv_2$. Subtract the results to get: $\int
\langle \sym \nabla \bfv_1 - \sym\nabla \bfv_2 : \nabla (\bfv_1 - \bfv_2)\rangle =
0$. Consequently:  $\int |\sym\nabla (\bfv_1-\bfv_2)|^2 =
0$ and thus $\bfv_1-\bfv_2$ must be a rigid motion. 

{\bf 2.} Note that $\bfv$ is a critical point (necessarily a minimizer) of the problem (E) if and only if:
\begin{equation}\label{var}
\int_\Omega \langle \sym\nabla \bfv : \nabla \mathbf{w}\rangle ~\mbox{d}x ~=~ 0
\qquad \mbox{for all} \quad \mathbf{w}\in W^{1,2}(\Omega,{\R}^d) \quad \mbox{with}
\quad \div \mathbf{w} =0. 
\end{equation}
Taking divergence free test functions which are compactly supported in
$\Omega$ and integrating by parts in (\ref{var}), it follows that $\div (\sym \nabla
\bfv) = \nabla p$ in the sense of distributions in $\Omega$, for some
$p\in L^2(\Omega, {\R})$. Here we use the convention that the divergence operator acts on rows of a
square matrix. This yields the first equation in (\ref{EL}). 
In addition, one has
\begin{equation}\label{var2}
\int_\Omega \langle \big(\sym\nabla \bfv - p I\big): \nabla \mathbf{w}\rangle ~\mbox{d}x = 0
\qquad \mbox{ for all } \quad \mathbf{w}\in W^{1,2}(\Omega,{\R}^d) \quad \mbox{with}
\quad \div \mathbf{w} =0. 
\end{equation}

Let now $\varphi\in\mathcal{C}^\infty_c(\partial\Omega,{\R}^d)$ satisfy:
\begin{equation}\label{condi}
\int_{\partial\Omega} \langle \varphi, \mathbf{n}\rangle =0.
\end{equation} 
Then there exists an divergence-free test function
$\mathbf{w}$ with trace $\mathbf{w}=\varphi$ on $\partial\Omega$. 
It is well known (see \cite{T}) that,
 since $(\sym\nabla \bfv - pI )$
 together with its divergence are square integrable in $\Omega$,
the normal trace $(\sym\nabla \bfv - p I )\mathbf{n}$ is well
defined on $\partial\Omega$. By (\ref{var2}) it thus follows
\begin{equation*}
0~ =~ \int_\Omega \langle \big(\sym\nabla \bfv - p I \big): \nabla \mathbf{w}\rangle
~\mbox{d}x~ =~ \int_{\partial\Omega} \langle \varphi, (\sym\nabla \bfv - p I )\mathbf{n}\rangle.
\end{equation*}
Since every tangential $\varphi$ obeys (\ref{condi}), it follows that
the tangential component of the normal stress vanishes:
$\big((\sym\nabla \bfv - p I )\mathbf{n}\big)_{tan}=0$. On the other hand,
the normal part satisfies
$$\langle (\sym\nabla \bfv - p I )\mathbf{n}, \mathbf{n}\rangle = const. \qquad
\mbox{ on } \partial\Omega.$$
Absorbing the constant in $p$, we obtain the boundary condition in (\ref{EL}).

{\bf 3.} To show (iii), let $\bar{\bfv}\in W^{1,2}(\Omega)$  be a solution
to $\div \bar{\bfv} =g(u)$, satisfying the bound (see \cite{T})
\begin{equation}\label{bar}
\|\bar{\bfv}\|_{W^{1,2}(\Omega)}~\leq ~C \|g(u)\|_{L^2(\Omega)}~\leq ~C \|u\|_{L^2(\Omega)}.
\end{equation}
Using $\mathbf{w}=\bfv-\bar{\bfv}$ as test function
in (\ref{var}), one obtains:
$$\int_\Omega |\sym\nabla \bfv|^2 ~=~ \int_\Omega\langle \sym\nabla \bfv : \nabla \bar
{\bfv}\rangle~ \leq ~\|\sym\nabla \bfv\|_{L^2(\Omega)} \|\nabla \bar{\bfv}\|_{L^2(\Omega)}. $$
In view of Korn's inequality and of (\ref{bar}), this yields the bound on the first term in (\ref{basic}). 
Since $\nabla p =
\div(\sym\nabla \bfv)$, we also obtain $\|p-\avint p\|_{L^2} \leq
C\|\nabla \bfv\|_{L^2(\Omega)}$ (see again \cite{T}). 
This completes the proof in view of $g$ being Lipschitz and $g(0)=0$.
\endproof

The next lemma states the uniform Schauder's estimates for the
classical solution of (\ref{EL}).

\begin{lemma} \label{lemS2}
Let $\Omega\subset \R^d$  be an open, bounded and
simply connected set with $\C^{2,\alpha}$ boundary. Let
$g\in\mathcal{C}^2(\R)$ be such that $g(0) = 0$ and $g'$, $g''$ are bounded.
Then, the boundary value problem (\ref{EL}) on $\Omega$
satisfies the ellipticity and the complementarity boundary conditions \cite{ADN}. 
Therefore its classical solution $(\bfv, p)$ 
satisfies the a-priori bound
\bel{S2}
\|\bfv\|_{\mathcal{C}^{2,\alpha}(\Omega)} + \|p\|_{\mathcal{C}^{1,\alpha}(\Omega)}~\leq~ C\big(\|g(u)\|
_{\mathcal{C}^{1,\alpha}(\Omega)} +
\|\bfv\|_{\mathcal{C}^{0,\alpha}(\Omega)} 
+\|p\|_{\mathcal{C}^{0,\alpha}(\Omega)}\big),
\end{equation}
where the constant $C$ depends only on $\Omega$.
Moreover, for every $u\in\mathcal{C}^{1,\alpha}(\Omega)$ the
minimization problem (\ref{3}) has a unique solution
$\bfv\in\C^{2,\alpha}(\Omega, \R^d)$
 normalized by the conditions
\bel{vnz}\avint_{\Omega} \bfv~\mathrm{d}x = 0\,,\qquad\qquad  
\skew \displaystyle{\avint_{\Omega} \nabla \bfv~\mathrm{d}x  =0.}\eeq  
This solution satisfies
\bel{S3}
\|\bfv\|_{\mathcal{C}^{2,\alpha}(\Omega)} ~\leq~ C \|g(u)\|_{\mathcal{C}^{1,\alpha}(\Omega)}.
\end{equation}
Further, for every constant $M>0$ and every domain $\Tilde\Omega$ for
which there exists a homeomorphism $\Lambda:\Omega\to{\Tilde\Omega}$ with
$\|\Lambda\|_{\C^{2,\alpha}(\Omega)},
\|\Lambda^{-1}\|_{\C^{2,\alpha}(\Tilde\Omega)}\leq M$, the
corresponding bound (\ref{S3}) 
is valid with a uniform constant $C$ that depends only on $M$ (in
addition to $\Omega$ and $\alpha$ that are given in the problem).
\end{lemma}
\v
{\bf Proof.} {\bf 1.}  We denote the right hand side function in (\ref{EL}):
\begin{equation}\label{aiut}
U=g\circ u
\end{equation}
and observe that $u\in\mathcal{C}^{1,\alpha}(\Omega)$ implies
$U\in\mathcal{C}^{1,\alpha}(\Omega)$ in view of the assumptions on $g$.

Let $(\bfv,p)\in W^{1,2}\times L^2$ be
the weak solution to (\ref{EL}) whose existence follows from Lemma
\ref{lem1}. To deduce that actually $\bfv\in W^{2,2}$ and $p\in
W^{1,2}$, one employs the usual difference quotients estimates (see
\cite{GT} for scalar elliptic problems and \cite{Lad} for systems with
Dirichlet boundary conditions), provided that the system is
elliptic and satisfies the complementarity conditions on the
boundary. We check these in the next steps below, for a slightly more
general system with nonconstant coefficients. Then, a repeated application of the
classical a-priori estimate due to Agmon, Douglis and Nirenberg
\cite{ADN} Theorem 9.3, combined with a Sobolev 
embedding estimate, yields: 
$$ \|\bfv\|_{W^{2,q}(\Omega)} + \|p\|_{W^{1,q}(\Omega)}~\leq~ C\big(\|U\|
_{W^{1, q}(\Omega)} + \|\bfv\|_{W^{1,q}(\Omega)} 
+\|p\|_{L^{q}(\Omega)}\big), $$
for every $2\leq q<\infty$, since
$U\in\mathcal{C}^{1,\alpha}(\Omega)$ 
implies $U\in W^{1,q}(\Omega)$. Consequently, by Morrey's embedding we have $(\bfv,p)\in
\C^{1,\gamma}\times \C^{0,\gamma}(%\ov
\Omega)$ for every $0<\gamma<1$. Applying the Schauder estimates \cite{ADN} Theorem 10.5,
we finally arrive at (\ref{S2}).

Let now $\Lambda$ and $M$ be as in the statement of the
lemma. Let $(\tilde{\bfv}, \tilde p)$ be 
the solution to (\ref{EL}) on a perturbed domain $\Tilde \Omega$,
for some right hand side $\tilde U \in \mathcal{C}^{1,\alpha}({\Tilde\Omega})$. Then the composition
$ (v^1,\ldots, v^d, p) = (\bfv, p) = (\tilde{\bfv}, \tilde p) \circ \Lambda\in
\mathcal{C}^{2,\alpha} \times \C^{1,\alpha}(\Omega)$ solves
the following boundary value problem for a system of $d+1$ equations:
\bel{pullback2}\left\{\bega{ll}
\displaystyle{\frac{1}{2}\Big\langle\nabla^2 v^i :
(\nabla\Lambda)^{-1}(\nabla\Lambda)^{-1, T}\Big\rangle + \frac{1}{2}\Big\langle
\sum_{k=1}^d(\nabla\Lambda)^{-1,T}(\nabla^2 v^k)
(\nabla\Lambda)^{-1}e_k, e_i\Big\rangle } & \\
\displaystyle{ \quad+
\Big\langle \nabla v^i, \Delta (\Lambda^{-1})\circ\Lambda\Big\rangle + 
\mbox{trace}\Big((\nabla \bfv) (\nabla\partial_i
(\Lambda^{-1})\circ\Lambda)\Big) - \Big\langle (\nabla \Lambda)^{-1, T}\nabla p,
  e_i\Big\rangle } & \\
\qquad\qquad\qquad\qquad\qquad\qquad
\qquad\qquad\qquad\qquad\qquad\qquad\qquad\qquad
\displaystyle{ = 0 }  & x\in \Omega, \cr
\displaystyle{ \Big\langle \nabla \bfv : (\nabla
  \Lambda)^{-1, T}\Big\rangle~ =~ \tilde U\circ\Lambda} &x\in \Omega, \vspace{2mm}\cr
\Big( \sym\big((\nabla \bfv) (\nabla\Lambda)^{-1}\big) -
p I \Big) (\nabla\Lambda)^{-1,T} \bfn  =0.  & x\in\partial \Omega.\enda\right.
\eeq
Note that, when $\Lambda=id$ is the identity map, 
the system (\ref{pullback2}) reduces to (\ref{EL}).

{\bf 2.} To show ellipticity and boundary complementarity of
(\ref{pullback2}), we use the standard notation in \cite{ADN}. 
The principal symbol is the square operator matrix $L_\Lambda$ of dimension $(d+1)\times (d+1)$,
 given in the block form below. Its coefficients are polynomials in the variables
$\xi = (\xi_1\ldots \xi_d)$, corresponding to differentiation in
directions $e_1\ldots e_d$ in $\Omega$:
\begin{equation*}
\begin{split} 
L_\Lambda(\xi) & = \left[\begin{array}{c|c} {\displaystyle{\frac{1}{2}
   \Big\langle \xi\otimes\xi : (\nabla\Lambda)^{-1}(\nabla\Lambda)^{-1, T}\Big\rangle I  +
   \frac{1}{2} (\nabla\Lambda)^{-1,T}(\xi\otimes \xi)(\nabla\Lambda)^{-1} } }\vspace{1mm} &
  -(\nabla\Lambda)^{-1,T}\xi \\ \hline 
\big((\nabla\Lambda)^{-1,T}\xi \big)^T & 0\end{array}\right] \\
& = L\big((\nabla\Lambda)^{-1,T}\xi\big),
\end{split}
\end{equation*}
where the $(d+1)\times (d+1)$ polynomial matrix $L = L_{id}$ is defined as:
\begin{equation}\label{LL}
L(\xi) = \left[\begin{array}{c|c} {{\frac{1}{2}|\xi|^2 I  +
    \frac{1}{2}\xi\otimes \xi }} \vspace{1mm} & -\xi \\ \hline 
\xi^T & 0\end{array}\right].
\end{equation}
The first $d$ rows in the matrix $L_\Lambda$ correspond to the equations in:
$\mbox{div} (\sym\nabla\tilde{\bfv} - pI)=0$; to
these rows we assign weights $s=0$. The last row corresponds to the
equation $\div \tilde{\bfv} = g(u)$; we assign to it the weight $s=-1$. The first $d$ columns in
$L_\Lambda$ correspond to the components of $\bfv$; to these columns we assign
weights $t=2$. The last column corresponds to $p$; we assign to it the weight $t=1$. 

In order to check the ellipticity of the operator $L_\Lambda$, we need to compute the determinant of
$L_\Lambda(\xi)$. The 
determinant of a block
matrix, where $D$ has dimension $1\times 1$, can be written as
$$\det \left[\begin{array}{c|c} A & B\\ \hline C &
    D\end{array}\right] = (D+1) \det A  - \det(A+B\otimes C).$$
Hence
$$\det L(\xi) = \det\big(\frac{1}{2}|\xi|^2 I  + \frac{1}{2}\xi\otimes \xi \big)
- \det\big(\frac{1}{2}|\xi|^2 I  - \frac{1}{2}\xi\otimes \xi \big).$$
Further, if $B$ is a square matrix of rank $1$, then $\det(A+B) =
\det A + \langle \mbox{cof } A : B\rangle $. Hence
$$\det (|\xi|^2 I  + \xi\otimes \xi) = |\xi|^{2d} + |\xi|^{2(d-1)}
\langle  I  : \xi\otimes \xi\rangle = 2 |\xi|^{2d} \quad \mbox{ and }
\quad \det (|\xi|^2 I - \xi\otimes \xi)=0.$$
Consequently, we obtain the ellipticity condition:
\begin{equation}\label{form22}
\det L_\Lambda(\xi) ~= ~\det L((\nabla\Lambda)^{-1,T}\xi)~ =~
\frac{1}{2^{d-1}}|(\nabla\Lambda)^{-1,T}\xi|^{2d} \neq 0 \qquad 
\forall\xi\neq 0.
\end{equation}

The supplementary condition on $L_\Lambda$ is also satisfied: for any pair of
linearly independent vectors $\xi, \bar\xi\in{\R}^d$ the
polynomial $\det L_\Lambda(\xi + \tau\bar\xi)$ in the complex variable $\tau$,
has exactly $d$ roots $\tau^+_\Lambda(\xi, \bar\xi)$ with positive
imaginary parts. The roots of $\det L(\xi+\tau\bar\xi)$ are all equal to
$$\tau^+(\xi, \bar\xi) = \frac{1}{|\bar\xi|^2} \big(-\langle \xi,
\bar\xi\rangle + i (|\xi|^2|\bar\xi|^2 - \langle \xi, \bar\xi\rangle^2)^{1/2}\big).$$
Finally, we find the adjoint of $L(\xi)$ by a direct calculation:
$$L^{adj}(\xi) = \big(\det L(\xi)\big) L(\xi)^{-1} =
\frac{|\xi|^{2d}}{2^{d-1}} \left[\begin{array}{c|c} {{\frac{2}{|\xi|^2} I  -
    \frac{2}{|\xi|^4}\xi\otimes \xi }} \vspace{1mm} & \frac{1}{|\xi|^2}\xi \\ \hline
\frac{1}{|\xi|^2}\xi^T & 1\end{array}\right].$$
Naturally, the following formulas correspond to the change of variable $\Lambda$:
\begin{equation*}
L_\Lambda^{adj}(\xi) = L^{adj}\big((\nabla\Lambda)^{-1,T}\xi\big),\qquad
\tau_\Lambda^+(\xi, \bar\xi) = \tau^+\big((\nabla\Lambda)^{-1,T}\xi, (\nabla\Lambda)^{-1,T}\bar\xi\big).
\end{equation*}

{\bf 3.} We now want to verify the complementing boundary condition at a
point $P\in\partial\Omega$ and relative to any tangent vector $\eta\in
T_P(\partial\Omega)$ perpendicular to the unit normal $\mathbf{n}$ to $\partial\Omega$ at
$P$. The boundary operator matrix $B_\Lambda$ in 
(\ref{pullback2}) is of
dimension $d\times (d+1)$. It has the block form as below, where we
assign to each row the same weight $r=-1$: 
\begin{equation*}
\begin{split}
B_\Lambda(\xi; \mathbf{n}) & = \left[\begin{array}{c|c}
    {\displaystyle{\frac{1}{2} \Big\langle (\nabla\Lambda)^{-1,T}\xi, (\nabla\Lambda)^{-1,T}\mathbf{n}\Big\rangle  I  +
    \frac{1}{2}(\nabla\Lambda)^{-1,T} (\xi\otimes \mathbf{n})
    (\nabla\Lambda)^{-1} }} & - \displaystyle{
  (\nabla\Lambda)^{-1,T}\mathbf{n} } \end{array}\right] \\
& = B\big((\nabla\Lambda)^{-1,T}\xi; (\nabla\Lambda)^{-1,T}\mathbf{n}\big),
\end{split}
\end{equation*}
and where the polynomial matrix $B= B_{id}$ is defined as:
$$ B(\xi; \bar\xi)= \left[\begin{array}{c|c}
    {\displaystyle{\frac{1}{2} \langle \xi, \bar\xi\rangle I  +
    \frac{1}{2}\xi\otimes \bar\xi }} & - \bar\xi \end{array}\right].$$
Compute the product
\begin{equation}\label{form}
\begin{split}
& D_\Lambda(\xi;\mathbf{n})  = B_\Lambda(\xi; \mathbf{n}) L^{adj}_\Lambda(\xi) = 
D\big(\nabla\Lambda)^{-1,T}\xi; (\nabla\Lambda)^{-1,T}\mathbf{n}\big),
\vspace{2mm} \\
& D(\xi; \bar\xi)  = \frac{|\xi|^{2d}}{2^{d-1}}\left[\begin{array}{c|c}
    {\displaystyle{\frac{\langle\xi, \bar\xi\rangle}{|\xi|^2} I - 2
        \frac{\langle\xi, \bar\xi\rangle}{|\xi|^4}\xi\otimes \xi +
        \frac{2}{|\xi|^2}\mbox{skew}(\xi\otimes \bar\xi)}}
&  \displaystyle{\frac{\langle\xi, \bar\xi\rangle}{|\xi|^2}\xi - \bar\xi}\end{array}\right].
\end{split}
\end{equation}

The complementing boundary condition requires that, for any nonzero
tangent vector $\eta\in T_P(\partial\Omega)$, the
$d\times (d+1)$ matrix $D_\Lambda(\tau\mathbf{n} + \eta; \mathbf{n})$,
whose entries are polynomials in the complex variable $\tau$, has 
rows which are linearly independent modulo the  polynomial
\begin{equation}\label{form3}
M^+(\tau) = \big(\tau - \tau^+_\Lambda(\eta,\mathbf{n})\big)^d =
\big(\tau - \tau^+(\zeta,N)\big)^d.
\end{equation}
We use here the notation 
\begin{equation}\label{defim}
N \doteq (\nabla\Lambda)^{-1,T}\mathbf{n}, \qquad \zeta \doteq (\nabla\Lambda)^{-1,T}\eta.
\end{equation}

We will now directly reduce all the entries of $D_\Lambda(\tau\mathbf{n} +
\eta; \mathbf{n})$ by $M^+$ and prove that the reduced matrix of
coefficients at $\tau^0$ has rank $d$. In view of (\ref{form}), we obtain
\begin{equation}\label{form1}
\begin{split}
& D_\Lambda(\tau\mathbf{n} + \eta; \mathbf{n})  = D(\tau N +
\zeta; N) = \frac{|\tau N+\zeta|^{2(d-2)}}{2^{d-1}}
\times \\ & \times  \left[\begin{array}{c}
  \displaystyle{|\tau N + \zeta|^2 \langle \tau N +
     \zeta, N\rangle I - 2 \langle \tau N + \zeta, N\rangle (\tau N + \zeta)^{\otimes 2}
+ |\tau N+\zeta|^2 \big( N\otimes\zeta - \zeta\otimes N\big)}\\ \hline 
\displaystyle{|\tau N+\zeta|^2 \langle \tau N +  \zeta, N\rangle (\tau N + \zeta)^T- |\tau N+\zeta|^4 N^T}
\end{array}\right]^T.
\end{split}
\end{equation}
Observe that the vectors $\eta$, $\mathbf{n}$ are perpendicular, whereas
$\zeta$ and $N$, in general, are not. However, $\langle
\zeta, N\rangle = \langle \eta,
(\nabla\Lambda)^{-1}(\nabla\Lambda)^{-1,T}\mathbf{n}\rangle$ and since
the metric tensor $(\nabla\Lambda)^{-1}(\nabla\Lambda)^{-1,T}$ is
uniformly positive definite on $\Omega$, it follows that
\begin{equation}\label{form2}
|\langle \zeta, N\rangle | \leq \alpha |\zeta| |N|,
\end{equation}
with a universal constant $\alpha\in (0,1)$ that depends only on $M$.

Denote $a=\big(|\zeta|^2 |N|^2 - \langle\zeta, N\rangle\big)^{1/2}$,
which is a positive number because of (\ref{form2}). Writing for
simplicity $\tau^+=\tau^+(\zeta, N)$, we obtain
\begin{equation}\label{form6}
\tau^+-\overline{\tau^+} = \frac{2ia}{|N|^2}, \qquad
\langle\tau^+N+\zeta, N\rangle=ia.
\end{equation}
It is also easy to check that:
\begin{equation*}
\begin{split}
&|\tau N+\zeta|^{2(d-1)} = (\tau - \tau^+)^{d-1} (\tau -
\overline{\tau^+})^{d-1}  \equiv (\tau - \tau^+)^{d-1} (\tau^+ -
\tau^+)^{d-1} \quad \mbox{mod } M^+ \\
& \qquad\qquad\qquad\qquad \qquad\qquad\qquad\qquad 
\qquad = (\tau - \tau^+)^{d-1} \Big(\frac{2ia}{|N|^2}\Big)^{d-1} \quad
\mbox{mod } M^+,\\
& \langle \tau N+\zeta, N\rangle I \equiv \langle \tau^+ N+\zeta,
N\rangle I \quad \mbox{mod } (\tau-\tau^+) = ia I \quad \mbox{mod } (\tau-\tau^+),\\
& \tau N + \zeta \equiv \tau^+ N +\zeta \quad \mbox{mod } (\tau-\tau^+).
\end{split}
\end{equation*}
Therefore, by (\ref{form6}) we get the reduction of the last column of $D_\Lambda$:
\begin{equation}\label{form7}
D_\Lambda(\tau\mathbf{n}+\eta;\mathbf{n})e_{d+1} \equiv
(\tau-\tau^+)^{d-1} \mathcal{Z}_{d+1} \quad \mbox{mod } M^+, 
\end{equation}
where
$$\mathcal{Z}_{d+1}= \Big(\frac{2ia}{|N|^2}\Big)^{d-1} (ia) (\tau^+ N + \zeta). $$
In the next step we shall reduce the entries of 
$D_\Lambda(\tau\mathbf{n}+\eta;\bfv)_{d\times d}$ by $M^+$.
\v
{\bf 4.} Arguing  as above, and observing that $\zeta\otimes N - N\otimes\zeta = (\tau^+ N +
\zeta)\otimes N - N\otimes  (\tau^+ N +\zeta)$, we obtain
\begin{equation*}
\begin{split}
|\tau N+\zeta|&^{2(d-1)}  \Big( \langle\tau N+\zeta, N\rangle I +
\zeta\otimes N - N\otimes \zeta\Big)\\ & \equiv (\tau - \tau^+)^{d-1}\Big(\frac{2ia}{|N|^2}\Big)^{d-1} 
\big( ia I +\zeta\otimes N - N\otimes \zeta\big) \quad \mbox{mod } M^+\\
& = (\tau - \tau^+)^{d-1}\Big[\Big(\frac{2ia}{|N|^2}\Big)^{d-1} 
\big( ia I - N\otimes (\tau^+ N + \zeta\Big)\big)  +
\frac{1}{ia}\mathcal{Z}_{d+1}\otimes N\Big]\quad \mbox{mod } M^+.
\end{split}
\end{equation*}
On the other hand:
\begin{equation*}
\begin{split}
|\tau& N+\zeta|^{2(d-2)} \langle\tau N+\zeta, N\rangle (\tau
N+\zeta)^{\otimes 2} \\ & \equiv (\tau - \tau^+)^{d-1}\Big(\frac{2ia}{|N|^2}\Big)^{d-2}(ia) 
\Big[\frac{|N|^2d}{2ia} (\tau^+ N+\zeta)^{\otimes 2} + N\otimes  (\tau^+
N+\zeta) +  (\tau^+ N+\zeta)\otimes N\Big] \\
&  \qquad \qquad + (\tau - \tau^+)^{d-2}\Big(\frac{2ia}{|N|^2}\Big)^{d-2}(ia) 
(\tau^+ N+\zeta)^{\otimes 2} \quad \mbox{mod } M^+\\
& =  (\tau - \tau^+)^{d-1}\Big[\Big(\frac{2ia}{|N|^2}\Big)^{d-2}(ia) 
N\otimes  (\tau^+ N+\zeta) +  \mathcal{Z}_{d+1}\otimes \Big(
\Big(\frac{|N|^2d}{2ia}\Big)^2 (\tau^+ N+\zeta) + \frac{|N|^2}{2ia} N\Big)\Big]\\
&  \qquad \qquad + (\tau - \tau^+)^{d-2}\frac{|N|^2}{2ia}\mathcal{Z}_{d+1}\otimes
(\tau^+ N+\zeta) \quad \mbox{mod } M^+.
\end{split}
\end{equation*}
Concluding, we obtain
\begin{equation}\label{form8}
\begin{split}
&D_\Lambda(\tau\mathbf{n}+\eta;\mathbf{n})_{d\times d} \equiv
\mathcal{Z}_{d\times d} \quad \mbox{mod } M^+, \qquad \mbox{where:}\vspace{1mm}\\
&\mathcal{Z}_{d\times d}=  (\tau - \tau^+)^{d-1}\Big(\frac{2ia}{|N|^2}\Big)^{d-1} 
\Big [ia I + \Big(\frac{|N|^2}{2} - 1\Big) N\otimes (\tau^N+\zeta)\Big]
\\ & \qquad\qquad\qquad + (\tau - \tau^+)^{d-1}\mathcal{Z}_{d+1}\otimes
\Big[\Big(\frac{|N|^2d}{2ia}\Big)^{2}  (\tau^+ N+\zeta) + \frac{|N|^2+2}{2ia}N\Big] \\
&  \qquad \qquad\qquad + (\tau - \tau^+)^{d-2} \frac{|N|^2}{2ia}
\mathcal{Z}_{d+1}\otimes (\tau^+ N +\zeta).
\end{split}
\end{equation}
Consider now the reduced polynomial matrix of dimension $d\times (d+1)$:
$$\mathcal{Z}(\tau;\eta,\mathbf{n}) = \left[\begin{array}{c|c}
  \mathcal{Z}_{d\times d} & (\tau-\tau^+)^{d-1}\mathcal{Z}_{d+1}\end{array}\right],$$
where $\mathcal{Z}_{d\times d}$ and $\mathcal{Z}_{d+1}$ are given in
(\ref{form7}), (\ref{form8}). The complementing boundary condition
states precisely that $\mathcal{Z}$ has maximal rank (equal $d$) over
the field of complex numbers $\CC$. To validate this statement, it suffices to
check that the complex-valued matrix $\mathcal{Z}(0;\eta, \mathbf{n})$ is of
maximal rank. By performing elementary column operations and using the
fact that $\tau^+\neq 0$, we observe that $\mathcal{Z}(0;\eta,
\mathbf{n})$ is similar to:
\begin{equation}\label{form9}
\mathcal{Z}'(0;\eta, \mathbf{n}) =
(-\tau^+)^{d-1}\Big(\frac{2ia}{|N|^2}\Big)^{d-1}
\left[\begin{array}{c|c} \displaystyle{ ia I + \Big(\frac{|N|^2}{2}-1\Big) N\otimes (\tau^+N + \zeta) }
&  \tau^+N+\zeta \end{array}\right].
\end{equation}
We then compute, using (\ref{form6}):
\begin{equation*}
\begin{split}
\det\Big[ I + \frac{1}{ia} \Big(\frac{|N|^2}{2}-1\Big)N\otimes
(\tau^+N+\zeta)\Big] & = 1 + \mbox{trace}\Big(\frac{1}{ia} \Big(\frac{|N|^2}{2}-1\Big)N\otimes
(\tau^+N+\zeta)\Big) \\ & = 1+ \frac{1}{ia} \Big(\frac{|N|^2}{2}-1\Big)
\langle \tau^+N+\zeta, N\rangle = \frac{|N|^2}{2}.
\end{split}
\end{equation*}
Moreover,
\begin{equation}\label{form10}
|\det (\mathcal{Z}')_{d\times d}| 
~= ~\Big|\frac{2\tau^+a}{|N|^2}\Big|^{d(d-1)} a^d  \frac{|N|^2}{2}
~\neq~ 0.
\end{equation}
This establishes the validity of the ellipticity and the boundary
complementarity conditions for the system (\ref{pullback2}),
and thus in particular for the system (\ref{EL}).
\v
{\bf 5.} By the previous step, we can apply Theorem 10.5 in 
\cite{ADN} and obtain the estimate
\begin{equation}\label{esti3}
\|\bfv\|_{\C^{2,\alpha}(\Omega)} + \|p\|_{\C^{1,\alpha}(\Omega)}
\leq C \big( \|g(\tilde u\circ \Lambda)\|_{\C^{1,\alpha}(\Omega)} +
\|\bfv\|_{\C^{0,\alpha}(\Omega)} + \|p\|_{\C^{0,\alpha}(\Omega)} \big),
\end{equation}
where the constant $C$ (in addition to its dependence
on $\Omega$ and
$\alpha$)  depends only on an upper bound for the following quantities: the
$\C^{1,\alpha}$ norms of the coefficients of the highest order terms in the equations in (\ref{pullback2});
the $\C^{0,\alpha}$ norms of the coefficients of the lower order terms; the uniform ellipticity
constant $\lambda_\Lambda$; and the inverse of the minor constant
$\kappa_\Lambda$ (which is denoted in \cite{ADN} by the symbol $\Delta$). It is clear that
the former two quantities depend only on $M$. We now prove that the
bounds on $\lambda_\Lambda$ and $(\kappa_\Lambda)^{-1}$ also depend only on $M$.

Indeed, $\lambda_\Lambda$ is defined in terms of the inequalities
$$\frac{1}{\lambda_\Lambda}|\xi|^{2d} \leq \det L_\Lambda(\xi) \leq \lambda_\Lambda|\xi|^{2d}.$$
By (\ref{form22}) we  can thus take $\lambda_\Lambda =
2^{d-1}(\|(\nabla\Lambda)^{-1}\|_{\C^0}^{2d}
+\|\nabla\Lambda\|_{\C^0}^{2d})\leq 2^{d} M^{2d}$, valid for every $x\in\ov\Omega$.

On the other hand, the minor constant $\kappa_\Lambda$ is defined as follows. For any boundary point
$P\in\partial\Omega$ and any tangent unit vector 
$\eta\in T_P(\partial\Omega)$ at $P$, we write
$$\big[\mathcal{Z}(\tau;\eta,\mathbf{n})\big]_{ij} =
\sum_{s=0}^{d-1}q_{ij}^s\tau^s\qquad \mbox{ for } i=1\ldots d, ~~ j=1\ldots d+1.$$
Construct the matrix $Q=[q_{ij}^s]$, having $d$ rows: $i=1\ldots d$, and
$(d+1)d$ columns: $j=1\ldots d+1$, $s=0\ldots d-1$. Under the
complementing boundary condition, the rank of $Q$ equals $d$. Hence, if
$Q^1\ldots Q^K$ denote all the $d$-dimensional square minors of $Q$, 
one has
$$\max_{l=1\ldots K} |\det Q^l|~>~0.$$
The minor constant $\kappa_\Lambda$ is precisely the infimum of these
quantities, over all boundary points $P$
and all tangent unit vectors $\eta$ as above. 
Clearly, $\kappa_\Lambda>0$ and
$$\kappa_\Lambda ~\geq ~
\inf_{P\in \partial\Omega,~ \eta\in T_P(\partial\Omega),
~|\eta|=1} \big|\det (\mathcal{Z}'(0;\eta,\mathbf{n}))_{d\times d}\big|.$$
By (\ref{form10}) and the formula for $\tau^+(\zeta, N)$, we obtain
\begin{equation}\label{form11}
\frac{1}{\kappa_\Lambda} \leq \sup_{P\in \partial\Omega,~ \eta\perp\mathbf{n},
~|\eta|=1}\Big(\frac{|N|^4}{2 a}\Big)^{d(d-1)}\frac{1}{a^d}\frac{2}{|N|^2}.
\end{equation}
Recalling (\ref{defim}) and observing that $a\geq
(1-\alpha)^{1/2}|\zeta||N|$ in view of (\ref{form2}),  we conclude
that the quantity on the right 
hand side of (\ref{form11}) is bounded from above in terms of a
(positive) power of $M$.
This completes the proof of (\ref{esti3}), valid with 
a constant $C$ that depends only on $M$.
\v
{\bf 6.} We now show that (\ref{esti3}) can be improved to
\begin{equation}\label{esti4}
\|\bfv\|_{\C^{2,\alpha}(\Omega)} + \|p\|_{\C^{1,\alpha}(\Omega)}
\leq  C \|\tilde U\circ \Lambda\|_{\C^{1,\alpha}(\Omega)},
\end{equation}
where the constant $C$ depends only on $M$,
provided that $(\bfv, p)$ are normalized according to
\begin{equation}\label{zaz}
\avint_{\Omega}|\det\nabla\Lambda|{\bfv}~\mbox{d}x~=~0,\qquad
\mbox{skew}\avint_{\Omega}|\det\nabla\Lambda|(\nabla{\bfv})(\nabla\Lambda)^{-1}~\mbox{d}x~=~0,
\qquad \avint_{\Omega}p|\det\nabla\Lambda|~\mbox{d}x~=~0.
\end{equation}
As in the proof of Lemma \ref{lemS1}, we argue by contradiction. Assume there are sequences
of diffeomorphisms $\Lambda_n$ such that $\|\Lambda_n\|_{\C^{2,\alpha}},
\|\Lambda_n^{-1}\|_{\C^{2,\alpha}}\leq M$, and of normalized solutions
$(\bfv_n, p_n)$ to (\ref{pullback2}) with some $\tilde U_n \in \C^{1,\alpha}(\Lambda_n(\Omega))$, such that
\begin{equation}\label{esti5}
\|\bfv_n\|_{\C^{2,\alpha}(\Omega)} +
\|p_n\|_{\C^{1,\alpha}(\Omega)} =1 \qquad \mbox{and} \qquad \| 
\tilde U_n\circ\Lambda_n\|_{\C^{1,\alpha}(\Omega)}\leq \frac{1}{n}.
\end{equation}
We extract converging subsequences:  $\Lambda_n\to\Lambda$,
$\bfv_n\to\bfv$, and $p_n\to p$, as $n\to\infty$, in appropriate H\"older spaces with a 
fixed exponent $\beta\in (0, \alpha)$. The above implies (\ref{zaz})
and, since $(\bfv, p)$
solves the problem (\ref{pullback2}) with $\tilde U = 0$, by the
uniqueness of weak solutions on $\Tilde \Omega=\Lambda(\Omega)$ 
stated in Lemma \ref{lem1} (i), we obtain that $\bfv=0$ and $p=0$. Consequently, both
$\|\bfv_n\|_{\C^{0,\alpha}}$ and $ \|p_n\|_{\C^{0,\alpha}}$
converge to $0$, and by (\ref{esti3}) we get a contradiction with the
first assumption in (\ref{esti5}). Hence (\ref{esti4}) is proved.

Finally, we have 
$$\|\tilde{\bfv}\|_{\C^{2,\alpha}(\Tilde\Omega)}~\leq ~C \|\bfv\|_{\C^{2,\alpha}(\Omega)}\,,
\qquad\qquad \|\tilde U\circ\Lambda\|_{\C^{1,\alpha}(\Omega)}\leq C
\|\tilde U\|_{\C^{1,\alpha}(\Tilde\Omega)}\,,$$
with a constant $C$ depending only on $M$. In view of
(\ref{esti4}) and recalling (\ref{aiut}), this completes  the proof of the estimate
(\ref{S3}), with a constant independent of  the domain $\Tilde\Omega$.
%We see that (\ref{S2}) follows from
%Theorem 9.3 in \cite{ADN}. Finally, the dependence of $C$ on $\Omega$
%and $M$ is a direct consequence of the change of variable $\Lambda$,
%as described on pages 72-73 of \cite{ADN}.
\endproof

\subsection{Step 3: The growth of the domain $\Omega$}

\begin{lemma} \label{domgrowth}
Let $\Omega\subset\R^d$ be an open, bounded, smooth and simply
connected set, satisfying the uniform inner and
outer sphere condition with radius $2\rho>0$. Let $\varphi:\Sigma\to
(-\frac{\rho}{2}, \frac{\rho}{2})$ be a $\C^{2,\alpha}$ map, defining
the set $\Omega^\varphi$ as in (\ref{omp})-(\ref{definorm}). Let
$\bfv\in\C^{2,\alpha}({\Omega^\varphi}, \R^d)$ and define the new set:
\bel{Omh}\Omega_\epsilon~\doteq~\bigl\{ x+ \epsilon \bfv(x);~x\in \Omega^\varphi\bigr\}.\eeq
Then, there exists $\epsilon_0>0$, depending only on the upper bounds
of $\|\varphi\|_{\C^{2,\alpha}}$ and $\|\bfv\|_{\C^{2,\alpha}(\Omega^\varphi)}$, such that for every 
$\epsilon 
< \epsilon_0$ the following holds. The set $\Omega_\epsilon $ is open
and it can be represented as $\Omega_\epsilon  = \Omega^\psi$ for some
$\psi:\Sigma\to \R$ satisfying the bound:
\bel{psin}
\|\psi\|_{\C^{2,\alpha}}~\leq~\|\varphi\|_{\C^{2,\alpha}} + C \epsilon 
\|\bfv\|_{\C^{2,\alpha}(\Omega^\varphi)}.\eeq
The constant $C$ above depends only on the upper bounds of
$\|\varphi\|_{\C^{2,\alpha}}$ and $\|\bfv\|_{\C^{2,\alpha}(\Omega^\varphi)}$.
\end{lemma}

{\bf Proof.} {\bf 1.} Let $L$ be the Lipschitz constant of $\bfv $ on
${\Omega^\varphi}$. Since by Lemma \ref{lem2.3} we have
$\Omega^\varphi = \Lambda(\Omega)$ for some $\C^{2, \alpha}$
homeomorphism satisfying $\|\nabla\Lambda\|_{\mathcal{C}^0}\leq M$, it
follows by integrating along a curve connecting $x$ and $y$ in
$\Omega^\varphi$ that $|\bfv (x) - \bfv (y)|\leq C_\Omega M\|\nabla
\bfv \|_{\mathcal{C}^0} |x-y|$, where $C_\Omega$ depends only on the geometry of $\Omega$. Thus:
\begin{equation}\label{Lipma}
\|\nabla \bfv \|_{\mathcal{C}^0} ~\leq ~L~\leq~ C\|\nabla \bfv
\|_{\mathcal{C}^0},
\end{equation}
where $C$ depends only on $\|\varphi\|_{\C^{2,\alpha}}$ (we always suppress the dependence on
the referential $\Omega$). 

Define $\epsilon_0 \doteq  \frac{1}{2L}$. Then, for every
$\epsilon <\epsilon_0$, the map $id+\epsilon  \bfv $ is a $\C^{2,\alpha}$
homeomorphism between the open sets $\Omega^\varphi$ and (the automatically
open image) $\Omega_\epsilon $. This is so because the gradient $I
+\epsilon \nabla \bfv$ is invertible, implying the local 
$\C^{2,\alpha}$ invertibility of the map, whereas the map itself is an injection,
since $x+\epsilon  \bfv (x) = y+\epsilon  \bfv (y)$ yields $x=y$ in view of:
$$|x-y| ~= ~\epsilon |\bfv (x) - \bfv (y)|~\leq ~\epsilon  L |x-y|~ \leq ~\frac{\epsilon }{\epsilon_0}|x-y|.$$
In particular, we  observe that  $\partial\Omega_\epsilon  = \{x+\epsilon  \bfv
(x); ~ x\in\partial\Omega^\varphi\}$.
\v
{\bf 2.} We now construct $\psi$ so that $\Omega_\epsilon  =
\Omega^\psi$. By covering the boundary $\Sigma$ with finitely many charts, it suffices to 
consider the case where
$$\Omega~=~\{ (x_1,  x')=(x_1,x_2,\ldots, x_d)\in\R^d;~x_1<0\},
\qquad \Omega^\varphi~=~\{ (x_1,  x');~x_1<\varphi( x')\}.$$
Given $\bfv =(v^1,v') =(v^1,v^2\ldots, v^d) $ and $\epsilon >0$ as above, $\psi$ is defined by the relation
\bel{psup}
\psi\Big( x' + \epsilon v'(\vp( x'),  x') \Big)~ =~ \varphi( x') + 
\epsilon  v^1\bigl(\varphi( x'),  x'\bigr).
\eeq
The existence of $\psi$ and the bound (\ref{psin}) now 
follow by the implicit function theorem.
\endproof

\subsection{Step 4: Updating the density $w$} 

Before we continue with the discrete time set-up, let us motivate the
implicit definition (\ref{wk1}) by the following natural observation
regarding the transport equation (\ref{5}).

\begin{lemma}\label{motivation}
Let $\{\Omega(t)\}_{t\in [0,T]}$ be a Lipschitz continuous family of
sets with $\C^{2,\alpha}$ boundaries, defined as in (\ref{6}) through
a Lipschitz vector field $\bfv:\D=\{(t,x); ~ t\in [0,T], ~
x\in\Omega(t)\}\to\R^d$, satisfying $\bfv(t, \cdot)\in
\C^{2,\alpha}(\Omega(t),\R^d)$ for every $t\in [0,T]$. Denote
$\{\Lambda^t:\Omega(0)\to \Omega(t)\}_{t\in [0,T]}$ the corresponding $1$-parameter family of
diffeomorphisms given by the ODE:
\begin{equation}\label{difL}
\frac{\mathrm{d}}{\mathrm{d}t} \Lambda^t(x) = \bfv(t, \Lambda^t(x)),
\qquad \Lambda^0 = id.
\end{equation}
Assume that $w\in \C^{0,\alpha}(\D,\R)$ is a nonnegative density
function that satisfies (\ref{5}) in the weak sense (see (\ref{wsol})
for the precise definition). Then:
\begin{equation}\label{fufu}
w(t, \Lambda^t(x)) = \frac{w(0,x)}{\det\nabla\Lambda^t(x)} \qquad \mbox{for
all } x\in\Omega(0), ~ t\in [0,T].
\end{equation}
\end{lemma}
{\bf Proof.}
We will prove (\ref{fufu}) under the assumption $w\in \C^1(\D)$. The
general case of lower regularity will follow by a standard
approximation argument. Observe that, by (\ref{5}),
\begin{equation*}
\begin{split}
\frac{\mbox{d}}{\mbox{d}t}w(t,\Lambda^t(x)) & ~= 
~w_t(t, \Lambda^t(x)) +
\big\langle \nabla w(t,\Lambda^t(x)), 
\frac{\mbox{d}}{\mbox{d}t} \Lambda^t(x)\big\rangle
~=~ \big( w_t + \langle \nabla w, \bfv\rangle\big)
 (t, \Lambda^t(x)) \\
& ~=~ \big( w_t + \div(w\bfv) - w\div \bfv\big) (t, \Lambda^t(x))~ =~
 - (w\div\bfv) (t, \Lambda^t(x)).
\end{split}
\end{equation*}
On the other hand, using the formula
\begin{equation}\label{nice}
\frac{\mbox{d}}{\mbox{d}t}\det F(t) = \det F(t) \mbox{trace}\big(F'(t) F(t)^{-1}\big),
\end{equation}
valid for any matrix function $t\mapsto F(t)\in \R^{d\times d}$, we obtain
\begin{equation}\label{nice2}
\begin{split}
\frac{\mbox{d}}{\mbox{d}t}\det\nabla \Lambda^t(x) & = \big(\det\nabla \Lambda^t(x)\big)
\mbox{trace}\Big(\big(\frac{\mbox{d}}{\mbox{d}t}\nabla
\Lambda^t(x)\big)(\nabla\Lambda^t(x))^{-1}\Big) \\ & = \big(\det\nabla \Lambda^t(x)\big)
\mbox{trace}\Big(\nabla \bfv (t, \Lambda^t(x)) \nabla\Lambda^t(x)(\nabla\Lambda^t(x))^{-1}\Big)
\\ & = \big(\det\nabla \Lambda^t(x)\big) \div \bfv(t,\Lambda^t(x)).
\end{split}
\end{equation}
Consequently:
$$\frac{\mbox{d}}{\mbox{d}t}\big(\ln w(t,\Lambda^t(x))\big) = -
\frac{\mbox{d}}{\mbox{d}t}\big(\ln \det\nabla\Lambda^t(x)\big) = 
\frac{\mbox{d}}{\mbox{d}t}\big(\ln \frac{1}{\det\nabla\Lambda^t(x)}\big), $$
which directly yields (\ref{fufu}).
\endproof
\v
\begin{lemma} \label{density}
In the same setting of Lemma \ref{domgrowth}, let $w\in \C^{0,\alpha}(\Omega^\vp)$ be a non-negative density 
and let $u\in\mathcal{C}^{2,\alpha}(\Omega^\varphi)$ be the solution of
(\ref{2}) on the set $\Omega^\varphi$.
Then, there exists $\epsilon_0>0$ such that for every $\epsilon < \epsilon_0$, 
a new density $w_\epsilon$ is well defined on the set
$\Omega_\epsilon$ in (\ref{Omh}) by setting implicitly:
\bel{wu}
w_\epsilon(x+\epsilon \bfv(x))~ \doteq~{w(x)\over \det(I+\epsilon \nabla\bfv(x))}\,.
\eeq
Moreover, $w_\epsilon\geq 0$ and the following estimate holds: 
\bel{wca}
\|w_\epsilon\|_{\C^{0,\alpha}(\Omega_\epsilon)}~\leq~
(1 + C\epsilon)\|w\|_{\C^{0,\alpha}(\Omega^\vp)}. \eeq
Both the threshold $\epsilon_0$ and the constant $C$ above depend only on the upper bounds of
$\|\varphi\|_{\C^{2,\alpha}}$ and $\|\bfv\|_{\C^{1}(\Omega^\varphi)}$.
%and $\|w\|_{\C^{0,\alpha}(\Omega^\vp)}$.
\end{lemma}

{\bf Proof.}
Let $L$ be the Lipschitz constant of $\bfv$ on $\Omega^\vp$.
As observed in the proof of Lemma \ref{domgrowth}, the map
$x\mapsto x+ \epsilon \bfv(x)$ is a $\mathcal{C}^{2,\alpha}$
homeomorphism between $\Omega^\vp$ and $\Omega_\epsilon$.
Hence both the numerator and denominator in (\ref{wu}) are well defined
$\mathcal{C}^{2,\alpha}$ functions on $\Omega_\epsilon$, 
for all $\epsilon < \epsilon_0$ as long as $\epsilon_0\leq \frac{1}{2L}$.
By (\ref{wu}) the function $w_\epsilon$ is well defined and  
non-negative, provided that $\epsilon<\epsilon_0$.

By (\ref{Lipma}), the choice of
$\epsilon_0$ depends only on the upper bounds of the
quantities $\|\varphi\|_{\C^{2,\alpha}}$ and $\|\bfv\|_{\C^{1}(\Omega^\varphi)}$. 
Writing $\det (I+\epsilon\nabla \bfv(x)) = 1
+\epsilon\,\mathcal{O}(\|\nabla \bfv\|_{\mathcal{C}^0} + \|\nabla \bfv\|_{\mathcal{C}^0}^d)$,
we also deduce
\begin{equation}\label{num0}
0 \leq w_\epsilon(x) \leq (1+C\epsilon) \| w\|_{\C^{0}(\Omega^\vp)},
\end{equation}
for $\epsilon<\epsilon_0$ and $C$ as indicated in the statement of the Lemma.

It remains to estimate the H\"older constant of $w_\epsilon$. Using
(\ref{S1}) and the fact that:
$$|(x+\epsilon \bfv(x)) - (y+\epsilon \bfv(y))| ~ \geq ~(1-\epsilon L) |x-y|,$$
we obtain
\begin{equation*}
\begin{split}
\big|w_\epsilon(x+\epsilon \bfv(x)) & - w_\epsilon (y+\epsilon \bfv(y))\big| \\
& \ds \leq ~\frac{|w(u)-w(y)|}{\det(I+\epsilon \nabla \bfv(x))}
+ w(y) \left|\frac{1}{\det(I+\epsilon \nabla \bfv(x))}-
\frac{1}{\det(I+\epsilon \nabla \bfv(y))} \right|\\
& \leq ~ \ds [\nabla w]_{\alpha}|x-y|^\alpha \bigl(1+C\epsilon \|\bfv\|_{\C^1}\bigr)
 +  \|w\|_{\C^0} C \epsilon \|\bfv\|_{\C^1} |x-y| \\
&\leq ~\ds \Big([\nabla w]_\alpha + C \epsilon
\|w\|_{\mathcal{C}^0}\Big) |x-y|^\alpha \\
&\ds \leq~ \Big([\nabla w]_\alpha + C \epsilon \|w\|_{\mathcal{C}^0}\Big) \frac{
|(x+\epsilon \bfv(x)) - (y+\epsilon \bfv(y))|^\alpha}{(1-\epsilon L)^\alpha} \\
&\ds \leq ~ \Big([\nabla w]_\alpha + C \epsilon
\|w\|_{\mathcal{C}^0} (1+L)\Big) |(x+\epsilon \bfv(x)) - (y+\epsilon \bfv(y))|^\alpha,
\end{split}
\end{equation*}
since $(1-\epsilon L)^{-\alpha} \leq (1+2\epsilon L)^\alpha \leq 1+2\epsilon L$.
In view of (\ref{num0}), this yields (\ref{wca}).
\endproof

\section{Continuous dependence on data}
 \label{sec:4}
\setcounter{equation}{0}

As proved in Lemma \ref{lemS1} and Lemma \ref{lemS2}, the regularity
estimates (\ref{S1}) and (\ref{S3}) hold with a constant $C$ which 
is uniformly valid  for a family of domains $\Omega$, obtained via
diffeomorphisms with uniformly controlled $\mathcal{C}^{2,\alpha}$ norms.
In this section we study in more detail how the  solutions  
$u, \bfv$ of (\ref{2}) and (\ref{EL}) change, under small perturbations of $\Omega$.  
   
\begin{lemma}\label{LTLemma}
Let $\Omega\subset\R^d$ be an open, bounded and simply connected set
with $\mathcal{C}^{2,\alpha}$ boundary. Let $w\in
\mathcal{C}^{0,\alpha}(\Omega)$ be a nonnegative function. Then there
exists $\epsilon_0>0$ such that the following holds. Consider a
homeomorphism $\Lambda:\Omega\to\Tilde\Omega
=\Lambda(\Omega)$, satisfying:
$\|\Lambda - id\|_{\mathcal{C}^{2,\alpha}(\Omega)} \leq \epsilon_0$
and define $\tilde w\in\mathcal{C}^{0,\alpha}(\Tilde\Omega)$ by
$$\tilde w(\Lambda(x)) ~=~ \frac{w(x)}{\det\Lambda(x)} 
\qquad \mbox{for all } ~x\in\Omega.$$
Let $u$ be the solution to (\ref{2}) and $\bfv$ be
the solution to the minimization problem (\ref{3}), normalized as in
(\ref{vnz}). Likewise, let $\tilde u$ and $\tilde \bfv$ be the
corresponding solutions of these problems on $\Tilde\Omega$. 
Assume that $g\in\mathcal{C}^3(\R)$ with $g(0)=0$ and $g'$, $g''$,
$g'''$ uniformly bounded. Then
\bel{4.7}
\| \tilde u\circ\Lambda - u\|_{\C^{2,\alpha}(\Omega)} ~\leq~C
\|\Lambda-id\|_{\C^{2,\alpha}(\Omega)} \|w\|_{\C^{0,\alpha}(\Omega)}.
\eeq
and
\bel{4.77}
\|\tilde\bfv\circ\Lambda - \bfv\|_{\C^{2,\alpha}(\Omega)}~\leq~C
\|\Lambda-id\|_{\C^{2,\alpha}(\Omega)} \|w\|_{\C^{0,\alpha}(\Omega)}
\big(1+ \|w\|_{\C^{0,\alpha}(\Omega)}^2\big).
\eeq
Both the threshold $\epsilon_0$ and the constant $C$ above depend only
on the domain $\Omega$, and they are uniform for a family of domains
that are homeomorphic with controlled $\C^{2,\alpha}$
norms (as in the statements of  Lemmas~\ref{lemS1} and \ref{lemS2}).
\end{lemma}
\v
{\bf Proof.}
{\bf 1.} We first observe that, choosing $\epsilon_0>0$ sufficiently small, the map $\Lambda$
has a $\C^{2,\alpha}$ inverse $\Lambda^{-1}$. In addition, 
$\tilde w\in \C^{2,\alpha}(\Tilde\Omega
)$  is well defined, nonnegative,
and satisfies
\begin{equation}\label{la1}
\|\tilde w\|_{\C^{0,\alpha}(\Tilde\Omega
)}~\leq~ C \|w\|_{\C^{0,\alpha}(\Omega)}.
\end{equation}
The existence and uniqueness of the corresponding 
solutions $u$ and $\tilde u$ follow from Lemma \ref{lemS1}.
We regard $u^\sharp=\tilde u \circ\Lambda$ as an approximate solution
of (\ref{2}), and estimate the error
quantities $e_1, e_2$ in
\begin{equation*}
\left\{\begin{array}{ll} \Delta(u^\sharp -u) - (u^\sharp - u) = e_1 &
    \qquad x\in\Omega\\
\langle\nabla (u^\sharp - u ), \mathbf{n}\rangle = e_2 & \qquad x\in\partial\Omega.
\end{array}\right.
\end{equation*}
By (\ref{S1}) and (\ref{la1}) we obtain
\begin{equation}\label{la2}
\|u^\sharp\|_{\C^{2,\alpha}(\Omega)}\leq C \|\tilde
w\|_{\C^{0,\alpha}(\Tilde\Omega)}\leq C \|w\|_{\C^{0,\alpha}(\Omega)}.
\end{equation}
On the other hand, $u^\sharp$ solves the boundary value problem
(\ref{pullback}), where $A(x) = \big((\nabla\Lambda)^T\nabla\Lambda\big)^{-1}(x)$.
An explicit calculation yields:
\begin{equation}\label{la3}
\|A- I \|_{\C^{1,\alpha}(\Omega)} +
\|\nabla^2(\Lambda^{-1})\circ\Lambda\|_{\C^{0,\alpha}(\Omega)} \leq  C
\|\Lambda-id\|_{\C^{2,\alpha}(\Omega)}.
\end{equation}
Subtracting the equality $$\Delta u^\sharp - u^\sharp~ =~ (\Delta
u^\sharp - u^\sharp) - (\langle \nabla^2 u^\sharp : A\rangle + \langle
\nabla u^\sharp, \Delta (\Lambda^{-1})\circ\Lambda\rangle - u) -
\tilde w\circ \Lambda$$ from $\Delta u - u = -w$, we obtain
$$e_1 ~= ~- \langle \nabla^2 u^\sharp : (A - I)\rangle - \langle
\nabla u^\sharp, \Delta (\Lambda^{-1})\circ\Lambda\rangle  -
(\tilde w\circ \Lambda - w).$$
Hence, by (\ref{la3}) and (\ref{la2}), we obtain the bound
\begin{equation*}
\|e_1\|_{\C^{0,\alpha}(\Omega)} \leq C \|u^\sharp\|_{\C^{2,\alpha}(\Omega)} \|\Lambda - id\|_{\C^{2,\alpha}(\Omega)}
+ \|w\big(1-\frac{1}{\det\nabla\Lambda}\big)\|_{\C^{0,\alpha}(\Omega)} \leq
\|\Lambda - id\|_{\C^{2,\alpha}(\Omega)} \|w\|_{\C^{0,\alpha}(\Omega)}.
\end{equation*}
Likewise, computing the difference between the boundary conditions of $u^\sharp$ and $u$, we obtain
$$e_2 ~= ~ \langle \nabla u^\sharp, \mathbf{n} \rangle
~=~ - \langle \nabla u^\sharp, (A-I)\mathbf{n} \rangle.$$
Therefore (\ref{la2}) and (\ref{la3}) imply
\begin{equation*}
\|e_2\|_{\C^{1,\alpha}(\Omega)} \leq C \|u^\sharp\|_{\C^{2,\alpha}(\Omega)} \|\Lambda - id\|_{\C^{1,\alpha}(\Omega)}
\leq C \|\Lambda - id\|_{\C^{2,\alpha}(\Omega)} \|w\|_{\C^{0,\alpha}(\Omega)}.
\end{equation*}
By Theorem 6.30 in \cite{GT} it now follows
$$\|u^\sharp - u\|_{\C^{2,\alpha}(\Omega)}\leq C\big( \|u^\sharp - u\|_{\C^{0,\alpha}(\Omega)}
+ \|\Lambda - id\|_{\C^{2,\alpha}(\Omega)} \|w\|_{\C^{0,\alpha}(\Omega)}\big),$$
and the usual argument by contradiction, as in the proof of Lemma
\ref{lemS1}, yields the required bound on $\|\tilde u\circ\Lambda -
u\|_{\C^{2,\alpha}(\Omega)}$ in (\ref{4.7}).
\v
{\bf 2.} In order to estimate $\|\tilde \bfv\circ\Lambda -
\bfv\|_{\C^{2,\alpha}(\Omega)}$, let
$(\tilde\bfv, \tilde p)$ and $(\bfv, p)$ be the
normalized solutions to (\ref{EL}) on the domains
$\Tilde\Omega$ and $\Omega$, respectively.
Call $\bfv^\sharp=\tilde\bfv\circ\Lambda$, $p^\sharp = \tilde
p\circ\Lambda$.   We regard $(\bfv^\sharp, p^\sharp)$ as an 
approximate solution to (\ref{EL}). Indeed, it satisfies
the boundary value problem
\begin{equation}\label{la4}
\left\{\begin{array}{ll}
\div\big(\sym\nabla (\bfv^\sharp - \bfv) - (p^\sharp - p)I\big) = e_3 &\qquad x\in\Omega\\
\div (\bfv^\sharp - \bfv) = e_4 &\qquad x\in\Omega\\
\big(\sym\nabla (\bfv^\sharp - \bfv) - (p^\sharp - p)I\big)\bfn = e_5 &\qquad x\in\partial\Omega,\\
\end{array}\right.
\end{equation}
with error terms $e_3,e_4,e_5$. As in the proof of Lemma \ref{lemS2}, Theorem 10.5 in \cite{ADN} yields
\begin{equation}\label{la5}
\begin{split}
\|\bfv^\sharp &- \bfv\|_{\C^{2,\alpha}(\Omega)} + \|p^\sharp -
p\|_{\C^{1,\alpha}(\Omega)} \\ & \leq C \big(\|\bfv^\sharp -
\bfv\|_{\C^{0,\alpha}(\Omega)} + \|p^\sharp -
p\|_{\C^{0,\alpha}(\Omega)} + \|e_3\|_{\C^{0,\alpha}(\Omega)} +
\|e_4\|_{\C^{1,\alpha}(\Omega)}  + \|e_5\|_{\C^{1,\alpha}(\Omega)} \big).
\end{split}
\end{equation}
We claim  that (\ref{la5}) can be replaced by
\begin{equation}\label{la6}
\begin{split}
\|\bfv^\sharp  - \bfv\|_{\C^{2,\alpha}(\Omega)} & + \|p^\sharp -
p\|_{\C^{1,\alpha}(\Omega)} \\ & \leq C \bigg(\Big|\avint_\Omega (\bfv^\sharp -
\bfv)~\mbox{d}x \Big| + \Big|\avint_\Omega \skew \nabla (\bfv^\sharp -
\bfv)~\mbox{d}x \Big| + \Big|\avint_\Omega (p^\sharp - p)~\mbox{d}x
\Big| \\ & \qquad \quad  + \|e_3\|_{\C^{0,\alpha}(\Omega)} +
\|e_4\|_{\C^{1,\alpha}(\Omega)}  + \|e_5\|_{\C^{1,\alpha}(\Omega)} \bigg).
\end{split}
\end{equation}
Otherwise, we could find a sequence $(\bfv^\sharp_n - \bfv_n, p^\sharp_n - p_n)$
solving (\ref{la4}) with corresponding right hand sides $e_3^n$,
$e_4^n$ and $e_5^n$, and such that the left hand side of (\ref{la6}) equals $1$
for every $n$, while the quantities in the right hand side 
converge to $0$, as $n\to\infty$. 
Fix $\beta\in (0, \alpha)$. Extracting a subsequence, we deduce
that $\bfv^\sharp_n - \bfv_n$ and $p^\sharp_n - p_n$ converge in
$\C^{2,\beta}(\Omega)$ and $\C^{1,\beta}(\Omega)$, respectively, to
some limiting fields $V$, $P$, that solve the homogeneous problem
(\ref{la4}). Moreover, all the averages: $\displaystyle{\avint_\Omega V~\mbox{d}x,
  \avint_\Omega P~\mbox{d}x, \avint_\Omega \skew
\nabla V~\mbox{d}x}$, equal $0$. By uniqueness, this implies
 $V=0$ and $P=0$. Hence $\|\bfv^\sharp_n - \bfv_n\|_{\C^{0,\alpha}(\Omega)}$ and $\|p^\sharp_n -
p_n\|_{\C^{0,\alpha}(\Omega)}$ converge to $0$. But, this contradicts the
uniform estimate (\ref{la5}), since the left hand side always equals $1$.
\v
{\bf 3.} We now compute the error quantities $e_3, e_4, e_5$  in (\ref{la4}). Since
$(\bfv^\sharp, p^\sharp)=(v^{\sharp 1},\ldots, v^{\sharp d}, p^\sharp)$
solve the system (\ref{pullback}) on $\Omega$, one has
\begin{equation*}
\begin{split}
e_3^i  = & -\frac{1}{2}\big\langle \nabla^2 v^{\sharp i} : (A-I)\big\rangle -
\frac{1}{2} \left\langle \sum_{k=1}^d\big[ (\nabla\Lambda)^{-1,T} (\nabla^2
v^{\sharp k}) (\nabla \Lambda)^{-1} - \nabla^2 v^{\sharp k}\big] e_k
\,,~ e_i\right\rangle
\\ & - \big\langle \nabla v^{\sharp i}, \Delta(\Lambda^{-1})\circ
\Lambda\big\rangle - \mbox{trace}\big((\nabla\bfv^\sharp)
(\nabla\partial_i(\Lambda^{-1})\circ \Lambda)\big) + \big\langle
\big((\nabla \Lambda)^{-1}- I\big)^T\nabla p^\sharp, e_i\big \rangle, \\
e_4  = & - \big\langle \nabla \bfv^{\sharp} : \big((\nabla\Lambda)^{-1}- I\big)^T\big\rangle 
+ g(u^\sharp) -g(u),\\
e_5 = & - \frac{1}{2}(\nabla\bfv^\sharp)(A-I)\bfn - \frac{1}{2} \big[ (\nabla\Lambda)^{-1}
(\nabla \bfv^{\sharp}) (\nabla\Lambda)^{-1}-
\nabla\bfv^\sharp\big]^T\bfn + p^\sharp \big((\nabla\Lambda)^{-1}- I\big)^T\bfn.
\end{split}
\end{equation*}
Using (\ref{la3}) and the obvious bound $\|(\nabla\Lambda)^{-1} -
I\|_{\C^{1,\alpha}(\Omega)} \leq C \|\Lambda - id\|_{\C^{2,\alpha}(\Omega)}$, we obtain
\begin{equation*}
\begin{split}
\|e_3&\|_{\C^{0,\alpha}(\Omega)}  + \|\langle \nabla \bfv^\sharp :
((\nabla\Lambda)^{-1} - I)^T\rangle\|_{\mathcal{C}^{1,\alpha}(\Omega)}
+ \|e_5\|_{\C^{1,\alpha}(\Omega)}  \\ & \leq
C \big(\|\bfv^\sharp\|_{\C^{2,\alpha}(\Omega)} + \|p^\sharp\|_{\C^{1,\alpha}(\Omega)}\big)
\|\Lambda - id\|_{\C^{2,\alpha}(\Omega)} \leq
C \|\Lambda - id\|_{\C^{2,\alpha}(\Omega)} \|g(\tilde u)\|_{\C^{1,\alpha}(\tilde\Omega)}\,.
\end{split}
\end{equation*}
Here we used (\ref{esti4}) in
\begin{equation*}
\|\bfv^\sharp\|_{\C^{2,\alpha}(\Omega)} + \|p^\sharp\|_{\C^{1,\alpha}(\Omega)} 
\leq C\big( \|\tilde \bfv\|_{\C^{2,\alpha}(\tilde\Omega)} + \|\tilde p\|_{\C^{1,\alpha}(\tilde\Omega)} \big)
\leq C \|g(\tilde u)\|_{\C^{1,\alpha}(\tilde\Omega)}.
\end{equation*}
Similarly, we check that
\begin{equation*}
\begin{split}
& \Big|\avint_\Omega (\bfv^\sharp -
\bfv)~\mbox{d}x \Big| + \Big|\avint_\Omega (p^\sharp - p)~\mbox{d}x\Big| = 
\Big|\avint_\Omega \bfv^\sharp (\det\nabla\Lambda -1) ~\mbox{d}x \Big| +
\Big|\avint_\Omega p^\sharp (\det\nabla\Lambda -1)~\mbox{d}x\Big| \\ &
\qquad\qquad  \qquad\leq~ C\big(
\|\bfv^\sharp\|_{\C^{0}(\Omega)} + \|p^\sharp\|_{\C^{0}(\Omega)}\big) 
 \|\Lambda - id\|_{\C^{1}(\Omega)} \leq  C\|\Lambda -
 id\|_{\C^{1}(\Omega)} \|g(\tilde u)\|_{\C^{1,\alpha}(\tilde\Omega)},\\
& \Big|\skew \avint_\Omega \nabla (\bfv^\sharp -
\bfv)~\mbox{d}x \Big| ~=~ \Big|\skew \avint_\Omega (\nabla \bfv^\sharp)
\big((\det\nabla \Lambda)(\nabla\Lambda)^{-1} - I\big)~\mbox{d}x \Big| \\ &
\qquad \qquad \qquad \leq~
C \|\bfv^\sharp\|_{\C^{1}(\Omega)} \|\Lambda - id\|_{\C^{1}(\Omega)} \leq  C\|\Lambda -
 id\|_{\C^{2,\alpha}(\Omega)} \|g(\tilde u)\|_{\C^{1,\alpha}(\tilde\Omega)}.
\end{split}
\end{equation*}

To bound the expression $\|g(\tilde u)\|_{\C^{1,\alpha}(\tilde\Omega)}$, we estimate
\begin{equation*}
\begin{split}
\big|\nabla (g\circ\tilde u)(x) - \nabla (g\circ\tilde u)(y)\big| & \leq
|g'(\tilde u(x)) - g'(\tilde u(y)) |\cdot |\nabla \tilde u(x)| +
|g'(\tilde u(y)|\cdot |\nabla \tilde u(x) - \nabla \tilde u(y)|\\ & \leq C
\|g''\|_{\mathcal{C}^0}\|\nabla \tilde u\|^2_{\mathcal{C}^0(\tilde\Omega)}
|x-y| + \|g'\|_{\mathcal{C}^0}\|\nabla \tilde u\|_{\mathcal{C}^{0, \alpha}(\tilde\Omega)} |x-y|^\alpha  
\end{split}
\end{equation*}
and thus, by (\ref{la1})
\begin{equation*}
\begin{split}
\|g\circ\tilde u\|_{\mathcal{C}^{1,\alpha}(\tilde\Omega)} & \leq
C \big( \|g''\|_{\mathcal{C}^0}\|\tilde u\|^2_{\mathcal{C}^1(\tilde\Omega)}
+ \|g'\|_{\mathcal{C}^0}\|\tilde u\|_{\mathcal{C}^{1,
    \alpha}(\tilde\Omega)} \big) \\ & \leq C \|\tilde w\|_{\mathcal{C}^{0,\alpha}(\tilde\Omega)}
\big(1+ \|\tilde w\|_{\mathcal{C}^{0,\alpha}(\tilde\Omega)}\big) 
\leq  C \|w\|_{\mathcal{C}^{0,\alpha}(\Omega)}
\big(1+ \|w\|_{\mathcal{C}^{0,\alpha}(\Omega)}\big)
\end{split}
\end{equation*}
\v
{\bf 4.} To conclude estimating the right hand side of (\ref{la6}), we
need to deal with the term $\|g(u^\sharp) - g(u)\|_{\C^{1,\alpha}(\Omega)}$. We have
$$ \|g(u^\sharp) - g(u)\|_{\C^0(\Omega)}~\leq~C \|g'\|_{\C^0} \, \|u^\sharp-u\|_{\C^0(\Omega)}\,,$$
and
$$\bega{rl}\bigl|\nabla (g\circ u^\sharp)(x) -
\nabla (g\circ u)(x)\bigr|&\leq~\Big| \bigl(g'(u^\sharp(x))-g'(u(x))\bigr)\,\nabla u^\sharp(x)\Big| + \Big|g'(u(x))\,\bigl(\nabla u^\sharp(x)-\nabla u(x)\bigr)\Big| \\[3mm]
&\le ~ \|g''\|_{\C^0}\|u^\sharp-u\|_{\C^0(\Omega)} \|\nabla u^\sharp\|_{\C^0(\Omega)}
+ \|g'\|_{\C^0} \|\nabla u^\sharp-\nabla u\|_{\C^0(\Omega)}\,.\enda $$
Moreover
\begin{equation*}
\begin{split}\bigl|\nabla (g&\circ u^\sharp)(x) -
\nabla (g\circ u)(x)-\nabla (g\circ u^\sharp)(y) +\nabla (g\circ
u)(y)\bigr|\\ & 
=~\big| g'(u^\sharp(x))\,\nabla u^\sharp(x)  -g'(u^\sharp(y))\,\nabla u^\sharp(y) 
- g'(u(x))\,\nabla u(x) +g'(u(y))\,\nabla u(y) \big|\\ &
 = ~\big|  \bigl(g'(u^\sharp(x)) - g'(u^\sharp(y)\bigr) \nabla u^\sharp(x)
+g'(u^\sharp (y)) \bigl( \nabla u^\sharp(x) -\nabla
u^\sharp(y)\bigr)\,\\ &
\qquad\qquad \qquad\qquad\qquad 
-  \bigl(g'(u(x)) - g'(u(y)\bigr) \nabla u(x) - g'(u (y)) \bigl(
\nabla u(x) -\nabla u(y)\bigr)\big|\\ & 
\leq ~\big|  \bigl(g'(u^\sharp(x)) - g'(u^\sharp(y)\bigr) \bigl(\nabla u^\sharp(x)-
\nabla u(x)\bigr)\big| \\ & \qquad\qquad\qquad \qquad\qquad + \big|
\bigl(g'(u^\sharp(x)) - g'(u^\sharp(y)\bigr) - \bigl(g'(u(x)) -
g'(u(y)\bigr) \big| \cdot |\nabla u(x)| \\ &
\qquad\qquad\qquad \qquad\qquad 
+\big|\bigl( g'(u^\sharp (y))- g'(u(y))\bigr)  \bigl( \nabla
u^\sharp(x) -\nabla u^\sharp(y)\bigr)\big| \\ &
\qquad \qquad\qquad \qquad\qquad 
+ | g'(u (y))|\cdot \big| \bigl( \nabla^\sharp u(x) -\nabla^\sharp u(y)\bigr) -
\bigl( \nabla u(x) -\nabla u(y)\bigr)\bigr) \big| \\ &
\leq~ C |x-y|\, \|g''\|_{\C^0} \|\nabla u^\sharp\|_{\C^0(\Omega)}
\|\nabla u^\sharp-\nabla u\|_{\C^0(\Omega)}\\ &
\quad + C \Big(\|g''\|_{\C^0} \|\nabla u^\sharp-\nabla u\|_{\C^0(\Omega)} \,|x-y|  \|\nabla u\|_{\C^0(\Omega)} 
+ \|g'''\|_{\C^0} \|u^\sharp-u\|_{\C^0(\Omega)} \|\nabla u\|^2_{\C^0(\Omega)}
\,|x-y|\Big) \\ &
\quad + \|g''\|_{\C^0} \|u^\sharp-u\|_{\C^0(\Omega)} \, \|u^\sharp\|_{\C^{1,\alpha}(\Omega)}
|x-y|^\alpha + \|g'\|_{\C^0} \|u^\sharp-u\|_{\C^{1,\alpha}(\Omega)} |x-y|^\alpha,
\end{split}
\end{equation*}
where the constant $C$ may depend on the geometry of $\Omega$.
We used here the following representation, valid for all $x, y$ such that $[x,y]\subset\Omega$
\begin{equation*}
\begin{split}
\bigl(g'(u^\sharp&(x)) - g'(u^\sharp(y)\bigr) - \bigl(g'(u(x)) -
g'(u(y)\bigr) \\ & 
=~\int_0^1 {d\over ds} \Big(g'\bigl(u^\sharp (sx+(1-s)y)\bigr)-
g'\bigl(u(sx+(1-s)y)\bigr)\Big)\,\mbox{d}s \\ & 
=~\int_0^1g''(u^\sharp) \langle \nabla u^\sharp- \nabla u , x-y\rangle
\, \mbox{d}s+\int_0^1\bigl( g''(u^\sharp) -g''(u)\bigr) \langle 
\nabla u, x-y\rangle \, \mbox{d}s\,.
\end{split}
\end{equation*}

Consequently, by (\ref{4.7}), (\ref{la1}) and the estimates in Lemma \ref{lemS1}, we get
\begin{equation*}
\begin{split}
\|g\circ u^\sharp - g\circ u\|_{\mathcal{C}^{1,\alpha}(\Omega)} & \leq
C \|u^\sharp - u\|_{\mathcal{C}^{1,\alpha}(\Omega)}
\big( 1+ \| u^\sharp\|_{\mathcal{C}^{1,\alpha}(\Omega)} + \| u\|_{\mathcal{C}^{1}(\Omega)} +
\| u\|_{\mathcal{C}^{1}(\Omega)} ^2\big) \\ & \leq 
C \|u^\sharp - u\|_{\mathcal{C}^{1,\alpha}(\Omega)}
\big( 1+ \| w\|_{\mathcal{C}^{0,\alpha}(\Omega)}^2\big) \\ & \leq 
C \|\Lambda - id\|_{\mathcal{C}^{2,\alpha}(\Omega)} \| w\|_{\mathcal{C}^{0,\alpha}(\Omega)}  
\big( 1+ \| w\|_{\mathcal{C}^{0,\alpha}(\Omega)}^2\big). 
\end{split}
\end{equation*}
In view of the bounds in Step 3, the proof of (\ref{4.77}) is done.
\endproof

\section{Local existence  of solutions to the growth problem}  \label{sec:5}
\setcounter{equation}{0}

By a solution to the growth problem (\ref{2v}-\ref{3}-\ref{5}-\ref{6})
on some time interval $[0,T]$, $T>0$, we mean: 
\begi 
\item A Lipschitz continuous family of sets $\{\Omega(t)\}_{t\in [0,T]}$ with $\C^{2,\alpha}$ boundaries, 
\item A Lipschitz continuous velocity field $\bfv(t,x)$ 
defined on the domain:
\begin{equation}\label{DDD}
\D=\{(t,x);~t\in [0,T],~x\in\Omega(t)\},
\end{equation}
with $\bfv(t,\cdot)\in \C^{2,\alpha}(\Omega(t), \R^d)$ for every $t\in [0,T]$,
\item A nonnegative, $\C^{0,\alpha}$ regular continuous density function $w= w(t,x)$ defined in $\D$, 
\endi
for which the following holds.
\begi
\item[(i)] For every $t\in [0,T]$, the set $\Omega(t)$ is  determined  by (\ref{6}),
\item[(ii)] The density $w$ provides a weak solution to  the transport
  equation (\ref{5}), namely
\bel{wsol}
\begin{split}
\int_{[0,T]\times \R^d} w \eta_t + w\langle\bfv, \nabla
\eta\rangle~\mbox{d}t\mbox{d}x & + \int_{\R^d} w_0(x)\eta(0,x)~\mbox{d}x
= 0 \\ &  \quad \mbox{ for all } \eta\in \C^{\infty}_c\big(\D\cap ([0, T)\times \R^d)\big),
\end{split}
\eeq
\item[(iii)] For every $t\in [0,T]$, the vector field 
$\bfv(t,\cdot)$ on $\Omega(t)$ is a minimizer of (\ref{3}), while
$u(t, \cdot)$ is the minimizer of (\ref{2v}) with $w=w(t,\cdot)$.
\endi

\begin{theorem}\label{thm1}
Assume that the initial domain $\Omega_0\subset\R^d$ is an open,
bounded, simply connected set with $\C^{2,\alpha}$ boundary $\Sigma_0$, for some $0<\alpha<1$.
Assume that $g$ satisfies (\ref{gprop}). Then, given an initial nonnegative density $w_0\in \C^{0,\alpha}(\Omega_0)$, 
the problem (\ref{2v}-\ref{3}-\ref{5}-\ref{6}) has a  solution
on some time interval $[0,T]$, with  $T>0$.
\end{theorem}
\v
{\bf Proof.}   
{\bf 1.} By the assumed regularity of $\Sigma_0$, the set 
$\Omega_0$ satisfies the uniform inner and outer sphere condition with a
radius $3\rho>0$.    We construct a new smooth, referential domain $\Omega$ and a function $\vp_0 =
\varphi\in\C^{2,\alpha}(\Sigma)$, so that the assertions of Lemma
\ref{lem2.4} hold with $\varepsilon_0 = \rho/3$. In particular, we have
$\Omega_0 = \Omega^{\vp_0}$. Introduce the constants
\bel{M02}
M_\varphi~\doteq~ 1 + \|\varphi_0\|_{\C^{2,\alpha}}, \qquad\qquad  
M_w~\doteq~1+\|w_0\|_{\C^{0,\alpha}(\Omega_0)}
\eeq
where the first norm refers to a $\rho$-neighborhood $V_\rho$ of $\Sigma$, as in (\ref{definorm}).

Fix a time step $0<\epsilon <\epsilon_0$, where $\epsilon_0>0$ is chosen small enough, as in
Lemma \ref{domgrowth} and Lemma \ref{density}, in connection
with the upper bounds $\|\varphi\|_{\mathcal{C}^{2,\alpha}}\leq M_\varphi$,
$\|w\|_{\mathcal{C}^{0,\alpha}(\Omega^\varphi)}\leq M_w$ and
$\|\bfv\|_{\mathcal{C}^{2,\alpha}(\Omega^\varphi)}\leq C_0 M_w(1+ M_w)$. 
The constant $C_0$ is such that $\|u\|_{\mathcal{C}^{2,\alpha}} \leq C_0 \|w\|_{\mathcal{C}^{0,\alpha}} $
and $\|\bfv\|_{\mathcal{C}^{2,\alpha}}\leq C_0\|w\|_{\mathcal{C}^{0,\alpha}} (1+ \|w\|_{\mathcal{C}^{0,\alpha}})$
according to (\ref{S1}) and (\ref{S3}), and it depends
only on $M_\varphi$ through Lemma \ref{lem2.3}. Consider the discrete times $t_k=k\epsilon$.   
For each $k=0,1,2,\ldots$, given the set $\Omega_k$ and the scalar 
nonnegative function $w_k\in \mathcal{C}^{0,\alpha}(\Omega_k)$, we follow steps 1--4 of Section \ref{sec:3} and 
construct a new density $w_{k+1}$ on the new set $\Omega_{k+1}$. As in (\ref{omp}), we use the
representation with an appropriate $\varphi_k\in \C^{2,\alpha}$:
\begin{equation*}
\Omega_k ~=~ \Omega^{\varphi_k}~ =~ \bigl\{ x\in\R^d;~\delta(x)< \varphi_k(\pi(x))\bigr\}.
\end{equation*}
We claim that, as long as  $t_k$ remains in a sufficiently small
interval $[0,T]$, the norms $\|w_k\|_{\C^{0,\alpha}(\Omega_k)}$
and $\|\vp_k\|_{\C^{2,\alpha}}$ satisfy a uniform bound,
independent of the time step $\epsilon>0$, namely
\bel{ibd}
\|\varphi_k\|_{\C^{2,\alpha}}~\leq~ M_\vp \quad \mbox{and} \quad
\|w_k\|_{\C^{0,\alpha}(\Omega_k)}~\leq~M_w.
\eeq
Indeed, by Lemmas \ref{lemS1}, \ref{lem1} and \ref{lemS2}, we see that the Schauder estimates yield
\begin{equation}\label{import}
\begin{split}
\|u_k\|_{\mathcal{C}^{2,\alpha}(\Omega_k)} & ~\leq ~C_0 \|w_k\|_{\mathcal{C}^{0,\alpha}(\Omega_k)}, \\
\|\bfv_k\|_{\mathcal{C}^{2,\alpha}(\Omega_k)} & ~\leq ~C_0
\|w_k\|_{\mathcal{C}^{0,\alpha}(\Omega_k)}\big(1+ \|w_k\|_{\mathcal{C}^{0,\alpha}(\Omega_k)}\big).
\end{split}
\end{equation}
In turn, by Lemma~\ref{domgrowth}, the new domain has the form 
$\Omega_{k+1} = \Omega^{\varphi_{k+1}}$, with
\bel{B1k}
\begin{split}
\|\varphi_{k+1}\|_{\C^{2,\alpha}}~& \leq~\|\varphi_k\|_{\C^{2,\alpha}} + C\epsilon
\|\bfv_k\|_{\C^{2,\alpha}(\Omega_k)} \\ &
\leq~\|\varphi_k\|_{\C^{2,\alpha}} + CC_0(1+M_w)\epsilon
\|w_k\|_{\C^{0,\alpha}(\Omega_k)} ~\doteq~\|\varphi_k\|_{\C^{2,\alpha}} + C_1\epsilon
\|w_k\|_{\C^{0,\alpha}(\Omega_k)},
\end{split}
\eeq
while by Lemma \ref{density} the density $w_{k+1}$ on $\Omega_{k+1}$ satisfies the estimate
\bel{B2k}
\|w_{k+1}\|_{\C^{0,\alpha}(\Omega_{k+1})}~\leq~
\|w_k\|_{\C^{0,\alpha}(\Omega_k)} + C_2 \epsilon\,  \|w_k\|_{\C^{0,\alpha}(\Omega_k)}.
\eeq
The constants $C_1, C_2$ remain uniformly bounded, as long as
$\vp_k, w_k$ satisfy (\ref{ibd}).  
Let now
$$T~\doteq~\min\left\{ {1\over {C_1M_w}}\,,\,  {1\over {C_2 M_w}}\right\}.$$
By (\ref{M02}), (\ref{B1k}), (\ref{B2k}), the bounds
(\ref{ibd}) are valid as long as $t_k\in [0,T]$, regardless of $\epsilon<\epsilon_0$.
\v
{\bf 2.} We write $\Omega^\epsilon(t_k) = \Omega_k$ and
$w^\epsilon(t_k,\cdot) = w_k$ at the  times $t_k= k\epsilon$ for
$k=0,1,2,\ldots, \lfloor \frac{T}{\epsilon}\rfloor +1$. 
The sets
$\Omega^\epsilon(t)$ and the functions
$w^\epsilon(t,\cdot)$ are then defined
 for all $t\in [0,T]$,
by linear interpolation.  More precisely, for  
$t\in [t_k, t_{k+1}]$
we define
\bel{OT}
\begin{split}
 \Omega^\epsilon(t) & ~\doteq~\{ x+ (t-t_k)\bfv_k(x)\,;~~x\in
 ~\Omega_k\},  \\  w^\epsilon(t,\,x+ (t-t_k)\bfv_k(x)) & ~
\doteq~{w_k(x) \over\det\bigl(I+ (t-t_k)\nabla\bfv_k(x)\bigr)}.
\end{split}
\eeq
Clearly, each $w^\epsilon$ is Lipschitz continuous in $t$.
We claim that $w^\epsilon$ are  
uniformly H\"older continuous in both variables $t$ and $x$.
Indeed, the uniform bounds on the norms $\|\bfv_k\|_{\C^{2,\alpha}(\Omega_k)}$ 
(see (\ref{import}) and (\ref{ibd}))
imply the uniform Lipschitz continuity of $\bfv_k$ in $x$, with a 
Lipschitz constant independent of the time step $\epsilon>0$:
\bel{Lipuv}|
\bfv_k(x)-\bfv_k(y)| ~\leq~L |x-y|\,.
\eeq
Given an initial point $x_0\in\Omega_0$, let $t\mapsto x(t, x_0)$
be the  characteristic of (\ref{OT}), starting at $x_0$; that is the polygonal line defined inductively by:
\begin{equation*} \label{achar}
x(0, x_0) = x_0 \quad \mbox{ and } \quad x(t, x_0) = x(t_k, x_0) + (t-t_k)\bfv_k (x(t_k, x_0))
\quad\mbox{for } t\in [t_k, t_{k+1}],
\end{equation*}
so that:
$$\Omega^\epsilon(t) = \Big\{x(t,x_0); ~ x_0\in\Omega_0\Big\}.$$
By (\ref{Lipuv}), it follows that for every $t_k=k\epsilon\in[0,T]$
and $x_0, \bar x_0\in\Omega_0$, we have: $(1-\epsilon L)^k 
|\bar x_0- x_0|\leq|x(t_k,\bar x_0) - x(t_k, x_0)|\leq (1+\epsilon
L)^k|\bar x_0- x_0|$. This yields:
\bel{Ld}
\begin{split}
e^{-2Lt}|\bar x_0 - x_0|~ & \leq~(1-\epsilon L)^{t/\epsilon} |\bar x_0-x_0|\\ & \leq~
|x(t,\bar x_0) - x(t, x_0)| ~ \leq~ (1+\epsilon L)^{t/\epsilon} |\bar x_0 - x_0|
~\leq~e^{Lt}|\bar x_0- x_0| \\ & \qquad\qquad  \qquad\qquad 
\qquad \leq~e^{LT}|\bar x_0-x_0| \qquad \mbox{for } t\in [0,T],
\end{split}
\eeq
where the lower bound holds for all $\epsilon>0$ small enough,  
while the upper bound holds for every $\epsilon$.
Using (\ref{nice}) and the definition (\ref{OT}),  
we compute the derivative of $w^\epsilon$ along a characteristic
$x(\cdot, x_0)$:
\begin{equation}\label{wso1}
\begin{split}
\frac{\mbox{d}}{\mbox{d}t} w^\epsilon &(t, x(t, x_0)) ~ =~
\frac{\mbox{d}}{\mbox{d}t} \left( \frac{w_k\big(x(t_k, x_0)\big)}{\det
  \big(I + (t-t_k)\nabla\bfv_k(x(t_k, x_0))\big)}\right) \\ &~ = ~
- w^\epsilon(t,x(t,x_0)) \mbox{trace}\Big( \nabla\bfv_k(x(t_k, x_0))
\big(I + (t-t_k)\nabla\bfv_k(x(t_k, x_0))\big)^{-1}\Big) \\ &
~=~ - w^\epsilon(t,x(t,x_0)) \div\bfv^\epsilon_{tr}(t, x(t, x_0)),
\end{split}
\end{equation}
where we trivially extend the definition of $\bfv_k$ at $t_k$ to
$\bfv^\epsilon_{tr}(t,\cdot)$ on $\Omega^\epsilon(t)$, for every $t\in
[0,T]$, by simply transporting its value along the characteristics:
$$ \bfv^\epsilon_{tr}\bigl(t,x+(t-t_k)\bfv_k(x)\bigr) =~ \bfv_k(x) \qquad
\mbox{for} \quad t\in[t_k, t_{k+1}).$$
Note that $\bfv^\epsilon_{tr}$ is not continuous (in time) at $t=t_k$.
However we still have the uniform bound on its spacial derivatives:
$\|\bfv^\epsilon_{tr} (t,\cdot)\|_{\mathcal{C}^{2,\alpha}(\Omega^\epsilon(t))}\leq
M_\bfv$, independent of $ \epsilon<\epsilon_0$ and valid for
all $t\in [0,T]$. The last equality in (\ref{wso1}) now follows from
the identity
$$\nabla\bfv^\epsilon_{tr} (t, x(t,x_0)) ~=~ \nabla\bfv_k(x(t_k, x_0))  \Big(I
+ (t-t_k)\nabla\bfv_k(x(t_k, x_0))\Big)^{-1}.$$ 
From (\ref{wso1}) we  obtain the representation formula
\begin{equation}\label{wso}
w^\epsilon (t, x(t, x_0))~ = ~\exp\left\{-\int_0^t \div\bfv^\epsilon_{tr}(s, x(s, x_0))~\mbox{d}s\right\} w_0(x_0).
\end{equation}
Therefore, for any $\tau_1 \leq \tau_2$ and $x_0, 
\bar x_0\in\Omega_0$, we have the estimate
\begin{equation}\label{gro}
\begin{split}
& \bigl|w^\epsilon (\tau_2, x(\tau_2, \bar x_0)) -  w^\epsilon (\tau_1, x(\tau_1,
x_0))\bigg|\\ &  \leq ~\Bigl| \exp \left\{ -\int_0^{\tau_2} \div \bfv^\epsilon_{tr} (s, x(s,
  \bar x_0))~\mbox{d}s \right\} -\exp \left\{ -\int_0^{\tau_1} \div \bfv^\epsilon_{tr} (s, x(s, x_0))~\mbox{d}s
\right\} \bigg| w_0(\bar x_0) \\ & 
\quad + \exp \left\{ -\int_0^{\tau_1} \div \bfv^\epsilon_{tr} (s, x(s, x_0))~\mbox{d}s\right\}
\bigl|w_0(\bar x_0)-w_0(x_0)\bigr|.
\end{split}
\end{equation}
By the uniform $\C^{2,\alpha}$ bound on
$\bfv^\epsilon_{tr}(t,\cdot)$ and by (\ref{Ld}), 
the first term in (\ref{gro}) satisfies
\begin{equation*}
\begin{split}
& C \Bigl|\int_{\tau_1}^{\tau_2} \div \bfv^\epsilon_{tr} (s, x(s,
  \bar x_0))~\mbox{d}s \Bigr| w_0(\bar x_0) + C \Bigl| \int_0^{\tau_1} \div
  \bfv^\epsilon_{tr} (s, x(s, \bar x_0)) - \div\bfv^\epsilon_{tr}(s, x(s,
  x_0))~\mbox{d}s\Bigr| w_0(\bar x_0) \\ & \qquad\qquad \leq 
C \|w_0\|_{\C^{0,\alpha}(\Omega_0)}\int_{\tau_1}^{\tau_2}
\|\bfv^\epsilon_{tr}(s,\cdot)\|_{\C^{2,\alpha}(\Omega^\epsilon(s))}~\mbox{d}s  
\\ & \qquad\qquad\quad + C \|w_0\|_{\C^{0,\alpha}(\Omega_0)} \int_{0}^{\tau_1}
\|\bfv^\epsilon_{tr}(s,\cdot)\|_{\C^{2,\alpha}(\Omega^\epsilon(s))} |x(s,
\bar x_0) - x(s, x_0)|~\mbox{d}s \\ & \qquad\qquad  
\leq C \left(\max_{t\in [0,
    T]}\|\bfv^\epsilon_{tr}(t,\cdot)\|_{\C^{2,\alpha}(\Omega^\epsilon(t))}\right)\|w_0\|_{\C^{0,\alpha}(\Omega_0)} 
\big(|\tau_1-\tau_2| + e^{LT}|\bar x_0 - x_0|\big).
\end{split}
\end{equation*}
Moreover,  the second term in (\ref{gro}) is bounded by $C \|w_0\|_{\C^{0,\alpha}(\Omega_0)}
|\bar x_0 -  x_0|^\alpha$. By (\ref{Ld}) we thus have
\begin{equation*}\label{gro2}
\begin{split}
\bigl|w^\epsilon (\tau_2, x(\tau_2, \bar x_0)) -  w^\epsilon (\tau_1, x(\tau_1,
x_0))\bigr| & ~\leq ~C \big(|\tau_1 - \tau_2|^\alpha + |\bar x_0 -
x_0|^\alpha\big) \\ & ~\leq~ C \big(|\tau_1 - \tau_2|^\alpha + |x(\tau_2, \bar x_0) -
x(\tau_1,x_0)|^\alpha\big),
\end{split}
\end{equation*}
where $C$ depends only on $M_w$, $M_\bfv$ and $T$, but it is
independent of $\epsilon$, as claimed.
\v
{\bf 3.} We now  examine the representation: $\Omega^\epsilon(t) =
\Omega^{\vp^\epsilon(t,\cdot)}$, where 
$\vp^\epsilon(t, \cdot)\in \mathcal{C}^{2,\alpha}(\Sigma)$ in view of Lemma \ref{lem2.4}. For
$t\in[t_k, t_{k+1}]$ we consider the homeomorphism 
$\Theta(t,\cdot): \Sigma\to\Sigma$, defined by 
$$\Theta(t,x)~\doteq~\pi\Big(x+\vp_k(x)\bfn(x) + (t-t_k) \bfv_k\big(x+\vp_k(x)\bfn(x)\big)\Big).$$
Observe that $\Theta(t,x)$ and $\Theta^{-1}(t,x)$ are uniformly Lipschitz continuous
in both $t$ and $x$. Since the map $\vp^\epsilon(t,\cdot):\Sigma\to\R$
can be implicitly defined by
$$ x+\vp_k(x)\bfn(x) +  (t-t_k) \bfv_k(x+\vp_k(x)\bfn(x))
~=~\Theta(t,x) + \vp^\ve(t,\Theta(t,x))\bfn(\Theta(t,x)),$$
it follows that $\vp^\epsilon$ is a Lipschitz  continuous function
of $(t,x)\in [0,T]\times \Sigma$, with a Lipschitz constant
independent of $\epsilon$.
\v
{\bf 4.} For every $t\in [0,T]$, we now define the velocity fields
$\bfv^\epsilon(t,\cdot)$ on $\Omega^\epsilon(t)$, by setting  
\bel{uwe2}
\begin{split}
\bfv^\epsilon\bigl(t,x+(t-&t_k)\bfv_k(x)\bigr)  ~\doteq ~\frac{t-t_k}{\epsilon} \bfv_{k+1}
(x+\epsilon\bfv_k(x)) + \big(1-\frac{t-t_k}{\epsilon} \big)\bfv_k(x)
\\ & =~ \frac{t-t_k}{\epsilon} \bfv_{k+1}
(x+\epsilon\bfv_k(x)) + \big(1-\frac{t-t_k}{\epsilon} \big)\bfv_{tr}^\epsilon(t,x+(t-t_k)\bfv_k(x)),
\end{split}
\eeq
whenever $t\in[t_k, t_{k+1})$ and $x\in\Omega_k$.
Notice that this provides an
 interpolation between the composition $\bfv_{k+1}\circ(id +
\epsilon\bfv_k)$ and $\bfv_k$, on $\Omega_k$. 
In view of (\ref{Ld}), it is clear that
$\|\bfv^\epsilon(t,\cdot)\|_{\mathcal{C}^{2,\alpha}(\Omega^\epsilon(t))}\leq
M_\bfv$, as before. 

We now claim that the vector fields $\bfv^\epsilon$ are uniformly
Lipschitz continuous in both variables $t$ and $x$. 
By Lemma \ref{LTLemma}, in view of (\ref{import}) and (\ref{ibd})
 we have the uniform bound
\begin{equation}\label{gl1}
\|\bfv_{k+1}\circ (id + \epsilon\bfv_k) - \bfv_k\|_{\C^{2,\alpha}(\Omega_k)} \leq 
C\epsilon \|\bfv_k\|_{\C^{2,\alpha}(\Omega_k)}
\|w_k\|_{\C^{0,\alpha}(\Omega_k)} \big(1+\|w_k\|_{\C^{0,\alpha}(\Omega_k)}^2\big) 
\leq C\epsilon\,.
\end{equation}
Observe that, for any $\tau_1\leq
\tau_2$ and $x_0, \bar x_0\in\Omega$,  one has
\begin{equation}\label{gl2}
\begin{split}
|\bfv^\epsilon&(\tau_2, x(\tau_2, \bar x_0)) - \bfv^\epsilon(\tau_1,
x(\tau_1, x_0))| \\ & \leq ~
|\bfv^\epsilon(\tau_2, x(\tau_2, \bar x_0)) - \bfv^\epsilon(\tau_1, x(\tau_1, \bar x_0))|
+ |\bfv^\epsilon(\tau_1, x(\tau_1, \bar x_0)) - \bfv^\epsilon(\tau_1,
x(\tau_1, x_0))|.
\end{split}
\end{equation}
To prove Lipschitz continuity in time, it is not restrictive
to assume that $\tau_1, \tau_2\in [t_k, t_{k+1}]$. Then, by
(\ref{uwe2}) and (\ref{gl1}) the first term on the right hand side of 
(\ref{gl2}) is bounded by
\begin{equation*}
\begin{split}
|\bfv^\epsilon(\tau_2, x(\tau_2, \bar x_0)) &- \bfv^\epsilon(\tau_1,
x(\tau_1, \bar x_0))| ~\\
&= ~\frac{\tau_2-\tau_1}{\epsilon} \left|\bfv_{k+1}
\big(x(t_k, \bar x_0) + \epsilon \bfv_k(x(t_k, \bar x_0))\big) -
\bfv_k(x(t_k, \bar x_0))\right| \\ & = ~\frac{\tau_2-\tau_1}{\epsilon} \left|\Big(\bfv_{k+1}\circ
(id + \epsilon \bfv_k) - \bfv_k\Big) (x(t_k, \bar x_0))\right| 
~ \leq~
C (\tau_2-\tau_1).
\end{split}
\end{equation*}
On the other hand,  in view of  (\ref{import}) and (\ref{ibd}), 
the second term in (\ref{gl2}) is bounded by
\begin{equation*}
\begin{split}
 |\bfv^\epsilon(\tau_1, x(\tau_1, \bar x_0)) &- \bfv^\epsilon(\tau_1,
x(\tau_1, \bar x_0))| \\ & \quad \leq~ \left|\bfv_{k+1}
\big(x(t_k, \bar x_0) + \epsilon \bfv_k(x(t_k, \bar x_0))\big) -
\bfv_{k+1}\big(x(t_k, \bar x_0) + \epsilon \bfv_k(x(t_k,
x_0))\big)\right| \\ &
\qquad\qquad
+ \left| \bfv_k(x(t_k, \bar x_0)) - \bfv_k(x(t_k,
x_0))\right| \\ & \quad \leq~ M_{\bfv}(2+\epsilon M_\bfv) |x(t_k, \bar x_0) - x(t_k, x_0)|.
\end{split}
\end{equation*}
Together, the above estimates yield a
Lipschitz bound on (\ref{gl2}):
\begin{equation*}
|\bfv^\epsilon(\tau_2, x(\tau_2, \bar x_0)) - \bfv^\epsilon(\tau_1,
x(\tau_1, x_0))| ~\leq ~C \Big(|\tau_1 - \tau_2| + |x(\tau_2,\bar
x_0) - x(\tau_1, x_0)|\Big).
\end{equation*}

In a similar way, we interpolate linearly along
characteristics and define the scalar function $u^\epsilon$ implicitly by setting
\begin{equation*}
u^\epsilon\bigl(t,x+(t-t_k)\bfv_k(x)\bigr) ~\doteq ~\frac{t-t_k}{\epsilon} u_{k+1}
(x+\epsilon\bfv_k(x)) + \big(1-\frac{t-t_k}{\epsilon} \big)u_k(x).
\end{equation*}
As in the previous case of $\bfv^\epsilon$, we conclude that the
norms $\|u^\epsilon(t,\cdot)\|_{\C^{2,\alpha}(\Omega^\epsilon(t))}\leq
M_u$ are uniformly bounded and that $u^\epsilon$ is 
uniformly Lipschitz continuous in both variables $t,x$.
\v

{\bf 5.} To avoid technicalities stemming from the fact that the
functions $w^\epsilon$, $u^\epsilon$,  $\bfv^\epsilon$ are defined on different
domains $\D^\epsilon = \bigl\{ (t,x);~ t\in [0, T], ~ x\in \Omega^\epsilon(t)\bigr\}$, we extend each of these maps
to the set $[0,T]\times B$, where $B\subset\R^d$ is a ball large enough to contain all
$\Omega^\epsilon(t)$.  By the analysis in previous steps, and the
appropriate uniform boundedness of $\vp^\epsilon, w^\epsilon$, $u^\epsilon$,
$\bfv^\epsilon$, the Ascoli-Arzel\`a compactness theorem, yields the
uniform convergence of (possibly subsequences, as $\epsilon_n\to 0$):
\begin{equation}\label{converg}
\begin{split}
 \vp^\epsilon\to \vp  \quad  \mbox{in } ~~ & \C^0([0,T]\times
 \Sigma, \R), \qquad\qquad  \bfv^\epsilon\to \bfv  \quad \mbox{in } ~~\C^0([0,T]\times B, \R^d)\\
& w^\epsilon\to w,\quad u^\epsilon\to u  \quad \mbox{in } ~~\C^0([0,T]\times B, \R)
\end{split}
\end{equation}
Defining $ \D = \bigl\{ (t,x);~ t\in [0, T], ~ x\in \Omega(t)\bigr\}$
as in (\ref{DDD}), where $\Omega(t) =\Omega^{\vp(t,\cdot)}$, we see
that the limit functions have the following properties:
\begin{itemize}
\item $\vp$ is Lipschitz continuous on $[0,T]\times \Sigma$ and
satisfies $\|\vp(t,\cdot)\|_{\C^{2,\alpha}}~\leq~M_\vp$  for all $ t\in [0, T]$, 
\item $w\in \C^{0,\alpha}(\D)$ is nonnegative and satisfies  $\|w(t,\cdot)\|_{\C^{0,\alpha}(\Omega(t))}\leq M_w$,
\item $u$ and $\bfv$ are Lipschitz continuous on $\D$ and satisfy the
  uniform bounds $\|u(t,\cdot)\|_{\C^{2,\alpha}(\Omega(t))}\leq M_u$,
$\|\bfv(t,\cdot)\|_{\C^{2,\alpha}(\Omega(t))}\leq M_\bfv$ for all $ t\in [0, T]$.
\end{itemize}
It remains to check the requirements (i)--(iii) in the definition of
solution to (\ref{2v}-\ref{3}-\ref{5}-\ref{6}).
To prove (i),  we first remark that %,  by the standard theory of ODEs, 
the uniform convergence of $\bfv^\epsilon$ in (\ref{converg}) implies 
the uniform convergence of $\bfv^\epsilon_{tr}$ to $\bfv$, because in
view of (\ref{gl1}) and (\ref{gl2}) we have:
$$\|\bfv^\epsilon(t,\cdot) - \bfv_{tr}^\epsilon(t,
\cdot)\|_{\C^0(\Omega^\epsilon(t))}~\leq ~\|\bfv_{k+1}\circ (id +
\epsilon\bfv_k) - \bfv_k\|_{\C^0(\Omega^\epsilon(t))}~ \leq ~
C\epsilon.$$
Consequently, the $\epsilon$-characteristics $t\mapsto x(t, x_0)$ that are trajectories
of the ODE
$$x'(t)~=~\bfv_{tr}^\epsilon(t,x(t)),\qquad\qquad x(0)=x_0\in\Omega_0,$$
converge, as $\epsilon\to 0$, to the corresponding trajectory of:
$$x'(t)~=~\bfv(t,x(t)),\qquad\qquad x(0)=x_0,$$
uniformly for $t\in [0,T]$. Note that $x(t)$ above is precisely given
by the diffeomorphisms in (\ref{difL}), with $x(t) = \Lambda^t(x_0)$.
Hence (\ref{6}) follows by (\ref{DDD}).

To prove (ii), we note that each $w^\epsilon$ is a weak solution of the linear transport equation:
$$w_t^\epsilon + \div(w^\epsilon\bfv_{tr}^\epsilon) = 0, \qquad \qquad w(0,\cdot) = w_0,$$
in view of (\ref{wso1}) and the identity
$${\frac{\mbox{d}}{\mbox{d}t}w^\epsilon(t, x(t,x_0)) 
~=~ w^\epsilon_t +
\Big\langle\nabla w^\epsilon, 
\frac{\mbox{d}}{\mbox{d}t}x(t, x_0)\Big\rangle~ =~
w_t^\epsilon + \langle \nabla w^\epsilon, 
\bfv^\epsilon_{tr}\rangle}.$$
Thanks to the uniform convergence in (\ref{converg}),  the limit density  $w$    provides a 
weak solution to the transport equation (\ref{5}), as expressed in (\ref{wsol}).

To prove (iii), we observe that $u(t,\cdot)$ is a minimizer of (\ref{2v}) if and only if
\bel{umi}\int_{\Omega(t)} 
 \langle\nabla u(t,x), \nabla \phi(x)\rangle + u(t,x) \phi(x) - w(t,x)\phi(x) ~\mbox{d}x~=~0,
\eeq
for every test function $\phi\in\C^\infty_c(\Omega(t))$.
Fix $t\in [0,T]$ and $\phi$ as above. By construction, there exists a
sequence of sets $\Omega^n=\Omega^{\vp^n}=\Omega^{\epsilon_n}(\tau_n)$, with 
$$\epsilon_n\to 0, \qquad
\tau_n=k_n \epsilon_n\to t\qquad\vp^n\to \vp(t, \cdot)\qquad \hbox{as} ~~n\to\infty.$$ 
Moreover, there exist functions 
$u^n= u^{\epsilon_n}(\tau_n,\cdot)$, $w^n= w^{\epsilon_n}(\tau_n, \cdot)$ on
$\Omega^n$, converging uniformly to $u(t,\cdot)$ 
and $w(t,\cdot)$ on every compact subset of $\Omega(t)$, such that
$$\int_{\Omega^n}  \langle\nabla u^n, \nabla \phi\rangle +
 u^n\phi - w^n\phi~\mbox{d}x~=~0.$$
Passing to the limit with $n\to\infty$ and recalling that $\nabla u^n$
converges to $\nabla u(t, \cdot)$, we get (\ref{umi}).

Likewise, there exists a sequence $\bfv^n=\bfv^{\epsilon_n}(\tau_n,
\cdot)$, converging uniformly to $\bfv(t, \cdot)$ on any compact
subset of $\Omega(t)$, and satisfying
\begin{equation*}
\int_{\Omega^n}\langle\bfv^n(x), \nabla \phi(x)\rangle -(g\circ u^n)(x)\phi(x)~\mbox{d}x~=~0,
\end{equation*}
for every test function $\phi$, since $\div \bfv^n=g(u^n)$ in
$\Omega^n$. Passing to the limit as $n\to\infty$, we obtain that
$\div \bfv(t,\cdot)= g(u(t,\cdot))$ holds in its equivalent weak sense:
$$ \int_{\Omega(t)}\langle\bfv(t,x),\nabla \phi(x)\rangle -g( u(t,x))\phi(x)~\mbox{d}x~=~0.$$
Finally, we show that for every $t\in [0,T]$, the vector field
$\bfv(t,\cdot)$ is a minimizer of (\ref{3}).  As in (\ref{var}), this is equivalent to
\begin{equation}\label{varr}
\int_{\Omega(t)} \langle \sym\nabla \bfv(t, x) : \nabla \mathbf{w}(x)\rangle ~\mbox{d}x ~=~ 0,
\end{equation}
for all divergence-free vector fields $\bfw\in \C^1(\Omega(t), \R^d)$.
Let $\bfw$ be such a vector field.  By construction, we have: $\int_{\Omega^n}
\langle \sym\nabla \bfv^n : \nabla \mathbf{w}\rangle ~\mbox{d}x = 0,$ 
whereas the uniform convergence $\nabla \bfv^n\to
\nabla \bfv(t,\cdot)$ implies (\ref{varr}).  This concludes the proof of the local existence.  
\endproof
\v
\begin{remark} 
(i) In our construction scheme,  the
discrete approximations $\bfv_k$ are normalized according to
(\ref{vkno}). As a consequence, the same properties  are valid for the limiting solution:  
\bel{vno}
\avint_{\Omega(t)} \bfv(t,x)~ \mbox{d}x~= ~0, \qquad  \skew
\avint_{\Omega(t)} \nabla  \bfv(t,x)~\mbox{d}x~=~0 \qquad \qquad \mbox{for }
t\in [0,T].
\eeq 

(ii) Calling $\ov T$ the maximal time of existence of solutions,
the proof of Theorem \ref{thm1} suggests that either $\ov T=+\infty$, or else as $t\to \ov T-$, 
one of the following possibilities occurs:
\begi
\item $\|w(t,\cdot)\|_{\C^{0,\alpha}(\Omega(t))}\to +\infty$,
\item The inner or the outer sphere condition fails, namely
%(see Figure \ref{f:sg22})  
$${Rad}(t)~=~\min\Big\{\inf_{x\in \partial\Omega(t)}
\,R_{in}(x),~~\inf_{x\in \partial\Omega(t)} \,R_{out}(x)\Big\} \to 0,$$
where $R_{in}(x)$ is the inner radius of curvature of  $\Omega(t)$ at
a boundary point $x$, and $R_{out}$ is the outer curvature radius.
\endi
\end{remark}

%Is there any obstruction to the use of Korn's inequality or regularity theory
%on a more general flat manifold, such as the one in Fig.~\ref{f:sg22}, right?

%\begin{figure}[htbp]
%\centering
 %\includegraphics[scale=0.4]{sg22.eps}
 %\includegraphics[scale=0.4]{FIG/sg22.eps}
 %   \caption{{\small Center: when the boundary loses regularity, 
 %   the evolution cannot be continued in time.   Right: as long as the coefficients of the Riemann
   % metric remain smooth, even if they lead to a self-intersection,  the evolution equations are well-posed.}}
%\label{f:sg22}
%\end{figure}

\v
\section{Uniqueness of the normalized solutions} \label{sec:6}
\setcounter{equation}{0}
It is straightforward to check that if the sets $\{\Omega(t)\}_{t\in [0,T]}$ and the functions
$(t,x)\mapsto w(t,x), \bfv(t,x)$ provide a solution to the problem
(\ref{2v}-\ref{3}-\ref{5}-\ref{6}), then 
infinitely many other solutions can be constructed by superimposing rigid motions:
\begin{equation*}
\begin{split}
&\Tilde\Omega(t)~=~\big\{R(t) x + \bfb(t);~x\in \Omega(t)\big\}, \\
\tilde w\bigl(t,\, R(t) x + \bfb(t)\bigr)  ~=&~w(t,x), \qquad 
\tilde \bfv\bigl(t,\, R(t) x + \bfb(t)\bigr) ~=~R(t) \bfv(t,x) +R'(t) x + \bfb'(t).
\end{split}
\end{equation*}
Here, $t\mapsto R(t)\in SO(d)$ and $t\mapsto \bfb(t)\in\R^d$ 
define a smooth path of rigid motions $t\mapsto R(t)x + \bfb(t)$ with
$R(0)~=~I$, $\bfb(0)=0.$  
The corresponding function $\tilde u$ is then implicitly 
defined by the identity 
$$ \tilde u\bigl(t,\, R(t) x + \bfb(t)\bigr)~  = ~u(t,x).$$

Note that the normalisation (\ref{vno}) for $\bfv$ implies that
$$\avint_{\Tilde\Omega(t)} \tilde\bfv(t,x)~\mbox{d}x~ =~
R'(t)\avint_{\Omega(t)}x~\mbox{d}x + \bfb'(t),\qquad \skew
\avint_{\Tilde\Omega (t)} \nabla\tilde \bfv(t,x) ~\mbox{d}x ~=~ R'(t)R(t)^T,$$
Therefore, (\ref{vno}) holds for $\tilde\bfv$ if and only if $R(t)=I$ and
$\bfb(t)=0$ for all $t$.

The next result shows that the normalized solution is unique.

\begin{figure}[htbp]
\centering
\includegraphics[scale=0.55]{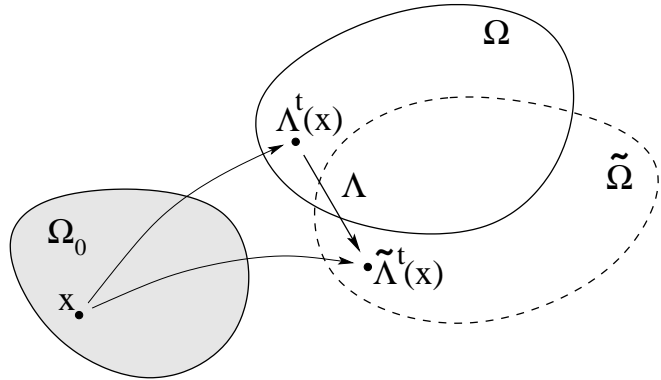}
    \caption{{\small The diffeomorphisms $\Lambda^t$ and $\tilde\Lambda^t$
      define the change of variable $\Lambda=\tilde\Lambda^t\circ (\Lambda^t)^{-1}$}}
\label{f:sg65}
\end{figure}

\begin{theorem}\label{thm2}
In the same setting as Theorem~\ref{thm1}, the problem 
(\ref{2v}-\ref{3}-\ref{5}-\ref{6}) has a unique solution which satisfies the additional identities
(\ref{vno}) for all $t\in [0,T]$.
\end{theorem}
\v
{\bf Proof.} Let $(\Omega, \bfv, w)$ and $(\Tilde\Omega, \tilde\bfv, \tilde w)$
be any two solutions, as defined in Section \ref{sec:5}, both
satisfying the normalization identities (\ref{vno}). For $t\in [0,T]$,
call $\Lambda^{t}:\Omega_0\to\Omega(t)$ and
$\Tilde\Lambda^{t}:\Omega_0\to\Tilde\Omega(t)$ the corresponding homeomorphisms
(see Figure \ref{f:sg65}) given by the ODEs  (\ref{difL}). We 
then have
\begin{equation}\label{z2}
\frac{\mbox{d}}{\mbox{d}t} \|\tilde\Lambda^t -
\Lambda^t\|_{\C^{2,\alpha}(\Omega_0)} ~\leq~
\|\tilde\bfv(t,\cdot)\circ\tilde\Lambda^t - \bfv(t,\cdot)\circ \Lambda^t\|_{\C^{2,\alpha}(\Omega_0)}.
\end{equation}
For a fixed $t\in [0,T]$, we shall apply Lemma \ref{LTLemma} to
the homeomorphism $\Lambda =
\tilde\Lambda^t\circ(\Lambda^t)^{-1}:\Omega(t) \to\Tilde\Omega
(t)$
and the nonnegative density $w(t,\cdot)\in\C^{0,\alpha}(\Omega(t))$.

The first assumption in Lemma \ref{LTLemma} holds for all sufficiently
small $t$, because
\begin{equation}\label{z1}
\|\Lambda- id\|_{\C^{2,\alpha}(\Omega(t))} ~= ~\|(\tilde\Lambda^t -
\Lambda^t)\circ (\Lambda^t)^{-1}\|_{\C^{2,\alpha}(\Omega(t))} ~\leq ~C
\|\tilde \Lambda^t - \Lambda^t\|_{\C^{2,\alpha}(\Omega_0)}~ \leq~\epsilon_0\,,
\end{equation}
because $\tilde\Lambda^0 = \Lambda^0= id$.
The second assumption follows by Lemma \ref{motivation}:
\begin{equation*}
\tilde w(t, \Lambda(x))~ =~
\frac{w_0\big((\Lambda^t)^{-1}(x)\big)}{\det\nabla\tilde\Lambda^t\big(
(\Lambda^t)^{-1}(x)\big)} ~= ~w(t,x) \frac{\det\nabla\Lambda^t\big(
(\Lambda^t)^{-1}(x)\big)}{\det\nabla\tilde\Lambda^t\big(
(\Lambda^t)^{-1}(x)\big)} ~= ~\frac{w(t,x)}{\det\nabla \Lambda(x)}.
\end{equation*}
Consequently, by (\ref{4.7}) we obtain
$$\|\tilde \bfv(t,\cdot)\circ\Lambda - \bfv(t,\cdot)\|_{\C^{2,\alpha}(\Omega(t))}~\leq~C \,
\| \Lambda - id\|_{\C^{2,\alpha}(\Omega(t))}.$$
Together with (\ref{z1}) this implies
\begin{equation*}
\begin{split}
\|\tilde \bfv(t,\cdot)\circ\tilde\Lambda^t - &\bfv(t,\cdot)\circ\Lambda^t\|_{\C^{2,\alpha}(\Omega_0)}
~ = ~\|\big(\tilde \bfv(t,\cdot)\circ\Lambda - \bfv(t,\cdot)\big)\circ\Lambda^t\|_{\C^{2,\alpha}(\Omega_0)}
\\ & \leq ~\|\tilde \bfv(t,\cdot)\circ\Lambda - \bfv(t,\cdot)\|_{\C^{2,\alpha}(\Omega(t))}~\leq ~
C \| \tilde\Lambda^t - \Lambda^t\|_{\C^{2,\alpha}(\Omega_0)},
\end{split}
\end{equation*}
for all times $t$ small enough, and with a uniform constant $C$.
Combining the above inequality with (\ref{z2}) we finally obtain
$${\mbox{d}\over \mbox{d}t} \| \tilde\Lambda^{t} - \Lambda^{t}
\|_{\C^{2,\alpha}(\Omega_0)}~\leq ~
C \| \tilde\Lambda^{t} - \Lambda^{t}
\|_{\C^{2,\alpha}(\Omega_0)}.$$
By Gronwall's inequality, this implies  that $\tilde\Lambda^t = \Lambda^t$
for all times  $t$ small enough.  In turn, this implies the 
equalities $\tilde w(t, \cdot) = w(t,\cdot)$ and
$\tilde u(t, \cdot) = u(t,\cdot)$. Likewise, $\tilde \bfv(t, \cdot)
= \bfv(t,\cdot)$, because of the normalization (\ref{vno}). Applying
the same argument on consecutive, sufficiently short time intervals,
we conclude that $(\Tilde\Omega, \tilde\bfv, \tilde w) = (\Omega, \bfv,
w)$ on the entire interval $t\in [0,T]$. \endproof

\section{Examples}\label{sec:75}
\setcounter{equation}{0}
We consider here two easy cases where the growth system can be solved explicitly.
\v
{\bf Example 1.} Assume that the volumetric growth rate
is proportional to the density of the morphogen, so that
$g(u)=au$ in (\ref{3}) with some $a>0$.  Then the volume of $\Omega(t)$ 
grows at a constant rate. Indeed, (\ref{6}) and (\ref{nice2}) give
\begin{equation*}
\begin{split}
{\mbox{d}\over \mbox{d}t} \hbox{vol}\,\Omega(t)~&=~
{\mbox{d}\over \mbox{d}t} \int_{\Omega_0} \det \nabla \Lambda^t(x)~\mbox{d}x
~=~ \int_{\Omega_0} \big(\det \nabla \Lambda^t(x)\big) \div\bfv(t,
\Lambda^t(x))~\mbox{d}x ~ \\ &=~ \int_{\Omega(t)} \div \bfv(t,x)~\mbox{d}x,
\end{split}
\end{equation*}
while from (\ref{2}) and since the conservation equation (\ref{5}) enjoys the solution formula
(\ref{fufu}), it follows that
\begin{equation*}
\int_{\Omega(t)} u(t,x)\, \mbox{d}x~=~\int_{\Omega(t)} (\Delta
u+w)(t,x)\, \mbox{d}x ~ =~\int_{\Omega(t)} w(t,x)\, \mbox{d}x~=~\int_{\Omega_0}
w_0(x)\, \mbox{d}x\,.
\end{equation*}
Concluding, the linear response function $g$ yields
\begin{equation}\label{dvol}
{\mbox{d}\over \mbox{d}t} \hbox{vol}\,\Omega(t)= a \int_{\Omega_0}
w_0(x)\, \mbox{d}x = a\kappa_0 \hbox{vol}\,\Omega_0 \qquad \mbox{ where } \quad
\kappa_0 \doteq\avint_{\Omega_0} w_0(x)\, \mbox{d}x\,.
\end{equation}
As a special case, assume that the initial domain $\Omega_0$ is a ball
centered at the origin with radius $r_0>0$, and the  initial density $w_0$ of signaling cells is radially symmetric.
By uniqueness (up to a rigid motion), the density $w(t,\cdot)$
remains then radially symmetric for all $t>0$, whereas the domain $\Omega(t)$ remains a ball
whose radius $r(t)$ may be determined from (\ref{dvol}), namely:
$r(t)^d=(1+\kappa_0 a t) r_0^d$. 

In particular, when $w_0(x) \equiv w_0>0$ is constant, then the quantities
\begin{equation}\label{przyklad1}
\begin{split}
& \Lambda^t(x) = (1+w_0 at)^{1/d}x, \qquad \Omega(t) = B\big(0, (1+w_0 at)^{1/d}\big),\\
& u(t,x) =w(t,x) = \frac{w_0}{1+w_0 a t},\\
& \bfv(t,x) = \frac{w_0 a}{d(1+w_0 a t)}x \quad \mbox{ and } \quad p(t,x)=\frac{w_0a}{d(1+w_0 a t)}
\end{split}
\end{equation}
provide the unique normalised solution to (\ref{2v}-\ref{3}-\ref{5}-\ref{6}).
\v
{\bf Example 2.} Next, assume that the growth rate $g:\R\to[0,\infty)$ is an arbitrary function satisfying
(\ref{gprop}), while the initial density $w_0$ of signaling cells is
again constant on an arbitrary domain $\Omega_0$ with center of mass at $0$,
so that: $\int_{\Omega_0} x~\mbox{d}x =0$.  In this case, for every $t\geq 0$
the density $w(t,x)= w(t)$ is spatially constant over the domain
$\Omega(t)$ and it satisfies the ODE
\bel{ODEw}
\dot w=- g(w) w\,,\qquad\qquad w(0)=w_0\,.\eeq
Indeed, generalizing (\ref{przyklad1}) we have that
\begin{equation*}
\begin{split}
& \Lambda^t(x) = \Big(\frac{w_0}{w(t)}\Big)^{1/d}x, \qquad \Omega(t) = \Big(\frac{w_0}{w(t)}\Big)^{1/d}\Omega_0,\\
& u(t,x) =w(t,x) = w(t),\\
& \bfv(t,x) = \frac{g(w(t))}{d}x \quad \mbox{ and } \quad p(t,x)=\frac{g(w(t))}{d}
\end{split}
\end{equation*}
solve (\ref{2v}-\ref{3}-\ref{5}-\ref{6}) together with (\ref{vno}). We further observe that setting:
$$w_{min}\doteq \max\{w\leq w_0; ~ g(w)=0\}\geq 0,$$
the solution to (\ref{ODEw}) satisfies $w(t)\to w_{min}$ as
$t\to\infty$. Consequently, if $w_{min}=0$ then $\Omega(t)$ becomes unbounded and
its volume approaches infinity. On the other hand, if $w_{min}>0$ then
$\Omega(t)$ increases to a finite limit $\Omega_{\infty} = \big(\frac{w_0}{w_{min}}\big)^{1/d}\Omega_0$.

\section{The Lagrangian formulation}\label{sec:7}
\setcounter{equation}{0}
In this section, we reformulate the coupled variational-transport
problem (\ref{2v}-\ref{3}-\ref{5}-\ref{6}) using the Lagrangian 
variable $\xi\in \Omega_0$ labeling points in the initial domain.

Let $\Lambda:[0,T]\times\Omega_0\to\R^d$ be the solution to the
problem in (\ref{6}), as in (\ref{difL}):
\begin{equation}\label{zero}
\frac{\mbox{d}}{\mbox{d}t}\Lambda(t,\xi) = \bfv(t, \Lambda(t,\xi)), \qquad \Lambda(0,\xi) = \xi.
\end{equation}
Define, for small $t\in[0,T]$, a flow of Riemann metrics 
$g:[0,T]\times \Omega_0\to {\R}^{d\times d}_{sym, >}$, by setting
\begin{equation}\label{trzy}
g(t,\xi) = \big((\nabla\Lambda)^T\nabla\Lambda\big) (t,\xi).
\end{equation} 
The Christoffel symbols of $g$ are given through:
$\partial_{ij}\Lambda = \sum_{m=1}^d\Gamma_{ij}^m\partial_m\Lambda$ or,
in vector notation: 
$$\Gamma_{ij}^\cdot = (\nabla \Lambda)^{-1}\partial_{ij}\Lambda\qquad
\forall i,j:1\ldots d.$$

We pull-back the solution quantities of the system
(\ref{2v}-\ref{3}-\ref{5}-\ref{6})  on  $\Omega_0$:
\begin{equation}\label{dwa}
 \tilde w(t,\xi) = w(t,\Lambda(t,\xi)), \quad \tilde u(t,\xi) =
u(t,\Lambda(t,\xi)), \quad \tilde{\bfv}(t,\xi) = \nabla\Lambda(t,\xi)^{-1} \bfv(t,\Lambda(t,\xi))
\end{equation}
and seek for their equivalent description (\ref{M1}-\ref{E1}-\ref{H1}-\ref{G1}) below.
There are some advantages in doing this:
\begi
\item A solution is a time-dependent field of $d\times d$ 
matrices $g= [g_{ij}]$ on
the fixed domain $ \Omega_0$.  
\item The transport equation (\ref{5}) has a trivial solution.
\item The non-uniqueness is automatically removed, since adding a
  rigid motion to the map $\xi\mapsto \Lambda(t,\xi)$ 
does not affect $g_{ij}$.
\item In Eulerian coordinates, the solution may cease to exist in
  finite time because different portions of the growing set may
  overlap. This issue does not arise when working in Lagrangian coordinates. 
\endi

On the other hand, while in Eulerian coordinates the elliptic 
equation
(\ref{2}) and the system (\ref{EL}) have 
constant coefficients, in Lagrangian coordinates
these coefficients depend on the metric itself. 
This makes the analysis considerably more difficult.

%\begin{figure}[htbp]
%\centering
%  \includegraphics[scale=0.5]{pg15.eps}
%    \caption{{\small Lagrangian and Eulerian coordinates.}}
%\label{f:pg15}
%\end{figure}

{\bf 1.} By Lemma \ref{motivation} and since $\det g = (\det\nabla\Lambda)^2$, we get:
\begin{equation}\tag{H1}\label{H1}
\tilde w(t,\xi) = \frac{ w_0(\xi)}{\sqrt{\det g(t,\xi)}}.
\end{equation}
To deal with (\ref{2v}), we observe equality of (the row) vectors in:
$\nabla u = (\nabla\tilde u)(\nabla \Lambda)^{-1}$, so that:
$$|\nabla u(t,\Lambda(t,\xi))|^2 = \big\langle (\nabla\tilde u)
(\nabla\Lambda)^{-1}(\nabla \Lambda)^{-1, T}, \nabla\tilde u\big\rangle =   \big\langle (\nabla\tilde u)
g^{-1}, \nabla\tilde u\big\rangle (t,\xi). $$
Changing the variables in (\ref{2v}) results in:
\begin{equation*}
\begin{split} 
J(u(t,\cdot)) & = \int_{\Omega_0}\Big(\frac{|\nabla u|^2}{2} + \frac{u^2}{2}
- wu\Big) (t, \Lambda(t,\xi)) \det\nabla \Lambda(t,\xi)~\mbox{d}\xi \\ & =
\int_{\Omega_0}\Big(\frac{1}{2} \big\langle (\nabla \tilde u) g^{-1} ,
\nabla \tilde u\big\rangle + \frac{1}{2}|\tilde u(t, \xi)|^2
- \tilde w\tilde u\Big) \sqrt{\det g(t,\xi)}~\mbox{d}\xi,
\end{split}
\end{equation*}
so that the minimization problem becomes:
\begin{equation}\tag{M1}\label{M1}
\mbox{minimize:} \qquad  \tilde J(t, \tilde u) =
\int_{\Omega_0}\Big(\frac{\big\langle (\nabla \tilde u) g^{-1} ,
\nabla \tilde u\Big\rangle}{2} + \frac{|\tilde u|^2}{2}
- \tilde w\tilde u\big)  \sqrt{\det g(t,\xi)}~\mbox{d}\xi.
\end{equation}

{\bf 2.} To rewrite (\ref{3}), differentiate the (column vector) equality
$\bfv(t,\Lambda(t,\xi))= (\nabla \Lambda)\tilde{\bfv}(t,\xi)$ in $\xi$:
\begin{equation}\label{pom5}
\begin{split}
& \nabla \bfv(t,\Lambda(t,\xi)) \\ & = (\nabla \Lambda)(\nabla\tilde
{\bfv})(\nabla\Lambda)^{-1}(t,\xi) + \Big[ (\partial_2\nabla\Lambda)\tilde{\bfv},
(\partial_1\nabla\Lambda)\tilde{\bfv}, \ldots,
(\partial_d\nabla\Lambda)\tilde{\bfv}\Big](\nabla\Lambda)^{-1}(t,\xi) \\ & = 
(\nabla \Lambda) \Bigg[ \nabla\tilde{\bfv}  + \Big[ (\nabla\Lambda)^{-1}(\partial_2\nabla\Lambda)\tilde{\bfv},
(\nabla\Lambda)^{-1} (\partial_1\nabla\Lambda)\tilde{\bfv}, \ldots,
(\nabla\Lambda)^{-1} (\partial_d\nabla\Lambda)\tilde{\bfv} \Big]\Bigg]
(\nabla\Lambda)^{-1}(t,\xi) \\ & 
= (\nabla \Lambda) (\tilde \nabla\tilde{\bfv}) (\nabla\Lambda)^{-1}(t,\xi), 
\end{split}
\end{equation}
where $\tilde \nabla\tilde{\bfv} = \{\tilde v^{i}_{,j}\}_{i,j=1\ldots d}$ is the covariant derivative of the
vector field $\tilde{\bfv} = \{\tilde v^i\}_{i=1\ldots d}$ with respect to
the metric $g$, in matrix notation given by:
$$\tilde\nabla\tilde{\bfv} = \nabla \tilde{\bfv}+
\Bigg[\Big[\Gamma_{11}^\cdot , \Gamma_{12}^\cdot , \ldots ,
\Gamma_{1d}^\cdot\Big] \tilde{\bfv}, \ldots , \Big[\Gamma_{j1}^\cdot , \Gamma_{j2}^\cdot , \ldots ,
\Gamma_{jd}^\cdot\Big]\tilde{\bfv}, \ldots, \Big[\Gamma_{d1}^\cdot , \Gamma_{d2}^\cdot , \ldots ,
\Gamma_{dd}^\cdot\Big] \tilde{\bfv} \Bigg],$$
so that $[\tilde\nabla\tilde{\bfv}]_{ij} = \tilde v^{i}_{,j} = \partial_j
\tilde v^i + \sum_{m=1}^d \Gamma_{jm}^i\tilde v^m$.
We thus obtain:
\begin{equation*}
\begin{split}
|\mbox{sym}\nabla \bfv|&^2 (t,\Lambda(t,\xi))   = ~  \frac{1}{4} \Big( \big\langle
(\nabla\Lambda)(\tilde\nabla\tilde{\bfv}) (\nabla\Lambda)^{-1} :
(\nabla\Lambda)(\tilde\nabla\tilde{\bfv}) (\nabla\Lambda)^{-1}\big\rangle  \\ &
\qquad \qquad \qquad \qquad + 2 \big\langle
(\nabla\Lambda)(\tilde\nabla\tilde{\bfv}) (\nabla\Lambda)^{-1} : 
(\nabla\Lambda)^{-1, T}(\tilde\nabla\tilde{\bfv})^T
(\nabla\Lambda)^{T}\big\rangle \\ & \qquad \qquad \qquad\qquad + 
\big\langle (\nabla\Lambda)^{-1, T}(\tilde\nabla\tilde{\bfv}) (\nabla\Lambda)^{T} :
(\nabla\Lambda)^{-1, T}(\tilde\nabla\tilde{\bfv})^T (\nabla\Lambda)^{T}\big\rangle\Big) \\ & = 
\frac{1}{2} \Big( \big\langle g(\tilde\nabla\tilde{\bfv}) g^{-1} :
\tilde\nabla\tilde{\bfv} \big\rangle  + 
\big\langle \tilde\nabla\tilde{\bfv} : (\tilde\nabla\tilde{\bfv})^T \big\rangle \Big) 
= \frac{1}{2} \Big( \big\langle g(\tilde\nabla\tilde{\bfv}) g^{-1} :
\tilde\nabla\tilde{\bfv} \big\rangle  + \mbox{trace} \big((\tilde\nabla\tilde{\bfv})^2\big)\Big).
\end{split}
\end{equation*}
Consequently, changing the variables in (\ref{3}) yields:
\begin{equation*}
\begin{split} 
E(\bfv(t,\cdot)) & = \frac{1}{2}\int_{\Omega_0}\big|\mbox{sym}\nabla\bfv (t,
\Lambda(t,\xi)) |^2 \det\nabla \Lambda(t,\xi)~\mbox{d}\xi \\ & = \frac{1}{4}
\int_{\Omega_0} \Big(\big\langle g(\tilde\nabla\tilde{\bfv}) g^{-1} :
\tilde\nabla\tilde{\bfv} \big\rangle  + \mbox{trace} \big((\tilde\nabla\tilde{\bfv})^2\big)\Big)(t,\xi)
 \sqrt{\det g(t,\xi)}~\mbox{d}\xi.
\end{split}
\end{equation*}
We further get:
$$\mbox{div } \bfv (t,\Lambda(t,\xi)) = \mbox{trace} \nabla \bfv
(t,\Lambda(t,\xi)) = \mbox{trace} \tilde\nabla\tilde{\bfv}(t,\xi) =
\widetilde{\mbox{div }} \tilde{\bfv} (t, \xi),$$
where the covariant divergence of the vector field $\tilde v$ is given by:
$$ \widetilde{\mbox{div }} \tilde{\bfv} =
\mbox{div}\nabla\tilde{\bfv} +
\sum_{k,i=1\ldots d} \Gamma_{ki}^k\tilde v^i = \div \nabla\tilde\bfv +
\langle \nabla\big( \ln\sqrt{\det g}\big), \tilde\bfv\rangle. $$
The minimization problem (\ref{3}) hence becomes:
\begin{equation}\tag{E1}\label{E1}
\begin{split}
\mbox{minimize:} \qquad & \tilde E(t,\tilde{\bfv}) 
= \frac{1}{4}\int_{\Omega_0} \Big(\big\langle g(\tilde\nabla\tilde{\bfv}) g^{-1} :
\tilde\nabla\tilde{\bfv} \big\rangle  + \mbox{trace} \big((\tilde\nabla\tilde{\bfv})^2\big)\Big)
 \sqrt{\det g(t,\xi)}~\mbox{d}\xi \\ & \mbox{ with }
\quad \widetilde{\mbox{div }}  \tilde{\bfv} = \tilde u. 
\end{split}
\end{equation}

We observe in passing that the integrand in (\ref{E1}) above depends only on the symmetric part of the
covariant derivative $\tilde\nabla\tilde{\bfv}_*$ of the
covariant tensor $\tilde{\bfv}_* =g\tilde{\bfv}$, 
carrying the resemblance to the original functional in (\ref{3}).
Indeed, since $\tilde \nabla\tilde{\bfv}_* = \tilde \nabla (g\tilde{\bfv}) = g \tilde
\nabla\tilde{\bfv}$, then $\tilde \nabla\tilde{\bfv}
=g^{-1} \tilde \nabla\tilde{\bfv}_*$, and:
\begin{equation*}
\begin{split}
\big\langle g(\tilde\nabla\tilde{\bfv}) & g^{-1} :
\tilde\nabla\tilde{\bfv} \big\rangle  + \mbox{trace}
\big((\tilde\nabla\tilde{\bfv})^2 \big) = 
\big\langle g^{-1}(\tilde\nabla\tilde{\bfv}_*) g^{-1} :
\tilde\nabla\tilde{\bfv}_* \big\rangle  + \mbox{trace}\big(
(g^{-1}\tilde\nabla\tilde{\bfv}_*)^2\big) \\ & 
= \big\langle g^{-1}(\tilde\nabla\tilde{\bfv}_*) g^{-1} :
\tilde\nabla\tilde{\bfv} \big\rangle  + \big\langle g^{-1}(\tilde\nabla\tilde{\bfv}_*) g^{-1} :
(\tilde\nabla\tilde{\bfv}_*)^T \big\rangle  \\ & =
2 \big\langle g^{-1}(\tilde\nabla\tilde{\bfv}_*) g^{-1} : \mbox{sym}\tilde\nabla\tilde{\bfv} \big\rangle =
2 \big\langle g^{-1}(\mbox{sym}\tilde\nabla\tilde{\bfv}_*)
g^{-1} : \mbox{sym}\tilde\nabla\tilde{\bfv} \big\rangle. 
\end{split}
\end{equation*}

{\bf 3.} The rule (\ref{6}) is being replaced by the equation for the evolution of the metric:
\begin{equation}\label{pom6}
\begin{split}
\frac{\mbox{d}}{\mbox{d}t} g(t, \xi) & = \frac{\mbox{d}}{\mbox{d}t}\big((\nabla \Lambda)^T\nabla\Lambda\big)(t, \Lambda(t,\xi))
 \\ & = \big(\nabla \bfv(t,\Lambda(t,\xi))\nabla \Lambda(t,\xi)\big)^T\nabla \Lambda +
 (\nabla\Lambda)^T \nabla \bfv(t,\Lambda(t,\xi))\nabla \Lambda(t,\xi) \\ &
= (\tilde\nabla\tilde{\bfv})^Tg + g (\tilde\nabla\tilde{\bfv}) = 2~\mbox{sym} \big(g
(\tilde\nabla\tilde{\bfv})\big)(t,\xi).
\end{split}
\end{equation}
We now conclude, by a direct calculation:
\begin{equation*}\tag{G1}\label{G1}
\begin{split}
\frac{\mbox{d}}{\mbox{d}t} g(t, \xi) = 2 ~\mbox{sym} \big(g
\nabla\tilde{\bfv}\big) + \sum_{i=1}^d (\partial_i g) \tilde v^i.
\end{split}
\end{equation*}
%\endproof

\section{Modeling the growth of a 2-dimensional surface in $\R^3$}
\label{sec:8}
\setcounter{equation}{0}
We now generalize the model (\ref{2v}-\ref{3}-\ref{5}-\ref{6}) to the case where,
instead of an open domain $\Omega(t)\subset\R^d$, the growing
set is  a codimension-one manifold $S(t)$. For simplicity, we assume
that $d=3$, so that $S(t)$ is a two-dimensional surface in $\R^3$.
\v
{\bf 1.}  Again, for each $t\in [0,T]$ we denote by $w(t,\cdot):S(t)\to\R$ a
nonnegative function representing the density of the signaling
cells in the tissue, whereas $u(t,\cdot):S(t)\to\R$ is the
concentration of produced morphogen. 
This function $u(t,\cdot)$ is defined to be the minimizer of
\begin{equation}\label{22v} \tag{M2}
\hbox{minimize:}\quad  J(u) = \int_{S(t)} \Big({|\nabla u|^2
\over 2} + {u^2\over 2} - wu\Big)~\mbox{d}\sigma(x),
\end{equation} 
or, equivalently, the solution to:
\begin{equation}\label{22}
\left\{\begin{array}{ll}
\Delta_{LB} u - u + w=0\qquad\qquad &x\in S(t)\cr
\langle\nabla u,\nu\rangle =0\qquad\qquad & x\in\partial S(t).
\end{array}\right.
\end{equation}
Here $\nu\in T_xS$ is the normal vector to the boundary $\partial S$,
and  $\Delta_{LB}u$ stands for the Laplace-Beltrami operator acting
on the scalar field $u$ on $S$. 

Consider a chart of $S$, so that $S=y(\omega)$ is parametrized by an
immersion $y:\omega\to\R^3$ for some open set $\omega\subset\R^2$. We
recall that the Laplace-Beltrami operator is given by
$$\Delta_{LB}u = \left[\displaystyle{\frac{1}{\sqrt{\det g}}}\sum_{i,j=1}^2
\partial_i\Big(\sqrt{\det g} ~g^{ij}\partial_j(u\circ y)\Big)\right]\circ y^{-1}.$$
On the domain $\omega$ of the chart, we denote by $[g_{ij}]_{i,j=1, 2}=(\nabla y)^T\nabla y $ 
the pull-back metric $g$ of the Euclidean metric
$I$ restricted to $S$, while its inverse is
denoted by $[g^{ij}]_{i,j=1, 2}= \bigl((\nabla y)^T\nabla y\bigr)^{-1}$.

\v
{\bf 2.} To determine the velocity $\bfv(t,\cdot):S(t)\to\R^3$, we
first derive the compressibility constraint expressing the fact that the infinitesimal
change of the surface area element due to the family of deformations
$\Lambda_\epsilon = id + \epsilon\bfv:S\to\R^3$ as $\epsilon\to 0$,
equals $u$.

Fix $t\in [0,T]$ and consider a flow of deformed surfaces $\epsilon\mapsto
\Lambda_\epsilon(S)$, starting from $S=S(t)$. 
For a given point $x\in S$, 
let $\{\tau_1(x), \tau_2(x)\}$ be an orthonormal  
basis of the tangent space $T_xS$. Calling $\bfn$ the unit normal vector to $S$, we compute
\begin{equation*}
\begin{split}
|\partial_{\tau_1}\Lambda_\epsilon \times \partial_{\tau_2}\Lambda_\epsilon|~
& =~ |(\tau_1 + \epsilon\partial_{\tau_1} \bfv)\times (\tau_2 + \epsilon \partial_{\tau_2}\bfv)| 
\\ & =~ |(\tau_1 \times \tau_2) + \epsilon(\partial_{\tau_1} \bfv \times \tau_2
-  \partial_{\tau_2}\bfv\times \tau_1) + \mathcal{O}(\epsilon^2)| \\ & =~
\Big(|\tau_1 \times \tau_2|^2 + 2\epsilon \big\langle \tau_1\times
\tau_2, \partial_{\tau_1} \bfv \times \tau_2 -  \partial_{\tau_2}\bfv\times
\tau_1\big\rangle + \mathcal{O}(\epsilon^2)\Big)^{1/2} \\ & = ~
  |\tau_1 \times \tau_2| \Big(1 + 2\epsilon \big\langle \frac{\tau_1\times
\tau_2}{|\tau_1\times \tau_2|^2}, \partial_{\tau_1} \bfv \times \tau_2 -  \partial_{\tau_2}\bfv\times
\tau_1\big\rangle + \mathcal{O}(\epsilon^2)\Big)^{1/2} \\ & =~
  |\tau_1 \times \tau_2| \Big(1 + \epsilon \big\langle \frac{\tau_1\times
\tau_2}{|\tau_1\times \tau_2|^2}, \partial_{\tau_1} \bfv \times \tau_2 -  \partial_{\tau_2}\bfv\times
\tau_1\big\rangle + \mathcal{O}(\epsilon^2)\Big) \\ & =~
|\tau_1 \times \tau_2| + \epsilon \big\langle \bfn, \partial_{\tau_1} \bfv \times \tau_2 -  \partial_{\tau_2}\bfv\times
\tau_1\big\rangle + \mathcal{O}(\epsilon^2).
\end{split}
\end{equation*}
By suitably choosing the orientation of $\bfn$, we can assume that
$\{\tau_1, \tau_2, \bfn\}$ is a positively oriented 
orthonormal basis of $\R^3$. Therefore 
\begin{equation*}
\begin{split}
\lim_{\epsilon\to 0}\frac{|\partial_{\tau_1}\Lambda_\epsilon
  \times \partial_{\tau_2}\Lambda_\epsilon| - |\tau_1\times \tau_2|}{\epsilon}~
& = ~\big\langle \bfn, \partial_{\tau_1} \bfv \times \tau_2 -  \partial_{\tau_2}\bfv\times
\tau_1\big\rangle \\ & =~ \big\langle \partial_{\tau_1} \bfv,  \tau_2 \times
\bfn\big\rangle -  \big\langle\partial_{\tau_2}\bfv,  \tau_1\times\bfn
\big\rangle\\ & =~  \big\langle \partial_{\tau_1} \bfv,  \tau_1 \big\rangle +  
\big\langle\partial_{\tau_2}\bfv,  \tau_2 \big\rangle.
\end{split}
\end{equation*}
We now decompose the vector field $\bfv=\bfv_{tan} + v_3\bfn$ into a tangential
component $\bfv_{tan}(x) \in T_xS$ and a normal component, given by a
scalar field $v_3:S\to\R$. Then
\begin{equation*}
\begin{split}
\big\langle \partial_{\tau_1} \bfv,  \tau_1 \big\rangle +  
\big\langle\partial_{\tau_2}\bfv,  \tau_2 \big\rangle & ~=~ \big\langle \partial_{\tau_1} \bfv_{tan},  \tau_1 \big\rangle +  
\big\langle\partial_{\tau_2}\bfv_{tan},  \tau_2 \big\rangle +
v_3\Big(\langle \partial_{\tau_1} \bfn,  \tau_1 \rangle +
\langle\partial_{\tau_2}\bfn,  \tau_2 \rangle\Big) \\ &~ = ~
\big\langle \partial_{\tau_1} \bfv_{tan},  \tau_1 \big\rangle +  
\big\langle\partial_{\tau_2}\bfv_{tan},  \tau_2 \big\rangle +
v_3\Big(\langle \Pi{\tau_1},  \tau_1 \rangle +
\langle \Pi{\tau_2},  \tau_2 \rangle\Big) \\ &~ = ~
\div \bfv_{tan} + v_3 \mbox{trace~} \Pi ~=~ \div \bfv_{tan} + 2H v_3,
\end{split}
\end{equation*}
where $\Pi =\nabla \bfn$ is the shape operator on $S$ and $H=
\frac{1}{2}\mbox{trace }\Pi$ is the mean curvature of $S$.
The constraint on $\bfv$ accounting for area growth can thus be written in the form
\begin{equation}\label{due.2}
\div \bfv_{tan} + 2H v_3~= ~u.
\end{equation}

To find an appropriate replacement of (\ref{3}) in the present setting, consider
the following model of elastic energy of deformations $\Lambda:S\to\R^3$ of $S$,
given by
$$I(\Lambda) ~=~ \int_S \mbox{dist}^2\big(\nabla \Lambda(x), O(2,3)\big)
~\mbox{d}\sigma(x).$$  
Here $O(2,3)=\{F\in\R^{3\times 2}; ~ F^TF=I\}$
represents gradients of deformations that preserve the metric on $S$. 
The integrand $\mbox{dist}^2(\cdot\,, O(2,3))$ may be replaced by
some other quadratic  function reflecting the material 
properties of the shell, provided it still satisfies the frame invariance 
and some other minimal regularity conditions.

Consider the expansion $\Lambda=id +\epsilon\bfv$. Then, in
analogy to the result in \cite{DNP}, we claim that the scaled functionals
$\epsilon^{-2}I$ $\Gamma$-converge as $\epsilon\to 0$ to the
following  elastic energy on $S$: 
\begin{equation}\label{due}
E(\bfv) ~= ~\frac{1}{2}\int_S |\sym\nabla \bfv_{tan} + v_3\Pi|^2~\mbox{d}\sigma(x).
\end{equation}
Among all velocity fields $\bfv$ which satisfy  (\ref{due.2}),
by the previous analysis we should thus choose one which minimizes 
(\ref{due.2}).   In the present setting, the 
constrained minimization  (\ref{3}) should be replaced by
\begin{equation}\label{33}\tag{E2}
\hbox{minimize:}\quad \int_S |\sym\nabla \bfv_{tan} + v_3\Pi|^2~\mbox{d}\sigma(x), \quad\qquad 
\hbox{subject to:}\quad \div \bfv_{tan} + 2H v_3= u.
\end{equation}

\v
{\bf 3.} The evolving surface $S(t)$ is now 
recovered as the %image of the Poincar\'e map at time $t$, 
set reached by trajectories of $\bfv$ starting in $S(0)$. Namely,
\begin{equation}\label{66}\tag{G2}
S(t)~=~
\bigg\{ \Lambda^t(x)\,;\quad \Lambda^t(0)=x\in S(0) ~\mbox{ and }~
 \frac{\mbox{d}}{\mbox{d}s}\Lambda^s(x)=\bfv(s, \Lambda^s(x)) ~~\forall s\in [0,t]\bigg\}.
\end{equation}
Again, the morphogen-producing cells are
transported along the flow, so that their density satisfies
\begin{equation}\label{55}\tag{H2}
w(t,\Lambda^t(x)) = \frac{w(0, x)}{\det \nabla \Lambda^t(x)}
\qquad  \mbox{for all } ~x\in S(0), ~ t\in [0,T],
\end{equation}
where $\det \nabla \Lambda^t(x)$ is the Jacobian of the linear map
$\nabla\Lambda^t(x):T_xS(0)\to T_{\Lambda^t(x)}S(t)$.

In conclusion, we propose (\ref{22v}-\ref{33}-\ref{66}-\ref{55})
as a model for thin shell/surface growth.   We 
leave the resulting system of PDEs as a topic for future study.

\begin{remark}
(i) In the flat case $S\subset\R^2$ and assuming the
in-plane evolution to the effect that $v_3=0$, the
constraint (\ref{due.2}) becomes: $\div \bfv = u$, which is precisely
the constraint in (\ref{3}).
In the general case, the infinitesimal change of area decouples into the
in-surface part $\mbox{div }\bfv_{tan}$, and 
% the change of area orthogonal to the surface given by 
$2H v_3$. Note that if $S$ is a minimal surface then all its variations (preserving the
boundary)  yield zero infinitesimal change of total area, so in
view of (\ref{due.2}) we get  $\int_S Hv_3 = 0$ for every $v_3$ vanishing on
$\partial S$. Thus $H\equiv 0$, as expected.

(ii) The problem (\ref{due.2}) is under-determined (one equation in three unknowns). 
% The natural regularity for, say $u\in L^2$, is $v_3\in L^2$
%and $\bfv_{tan}\in W^{1,2}$. 
Representing $\bfv_{tan}= \nabla \psi$ as the gradient of a scalar field $\psi$ on $S$, the equation (\ref{due.2}) 
can be replaced by the Laplace-Beltrami equation
$$\Delta_{LB} \psi ~= ~u - 2Hv_3\,.$$

(iii) The energy functional $E(\bfv)$ in (\ref{due}) measures
stretching, i.e. the change in metric on $S$ after the deformation to
$\Lambda_\epsilon(S)$, of order $\epsilon$. This functional can be augmented by
adding the bending term at a higher order:
\begin{equation}\label{due2}
\bar{E}(\bfv) = \frac{1}{2}\int_S |\sym\nabla \bfv_{tan} +
v_3\Pi|^2~\mbox{d}\sigma(x) + \frac{\mu}{24}
\int_S  |\big(\nabla ((\nabla\bfv)\bfn) - (\nabla\bfv)\Pi\big)_{tan}  
|^2~\mbox{d}\sigma(x).
\end{equation}
The integrand in the second term above measures the difference of
order $\epsilon$ between the shape operator $\Pi$  
on $S$ and the shape operator $\Pi_\epsilon$ of $\Lambda_\epsilon(S) =
id +\epsilon\bfv$. Alternatively, the tensor under this integral
represents the linear map: $T_xS\ni\tau\mapsto
\big(\partial_\tau(\nabla\bfv)\big)\bfn\in T_xS$. The presence of a
bending term introduces a regularizing effect, while 
the prefactor $\frac{\mu}{24}$, which is a fixed  small ``viscosity'' parameter,
guarantees that bending contributes at a higher order than stretching.

Let us also mention that a potentially relevant to the problem at hand
discussion of the $2$-dimensional models of elastic shells and their 
relation to the $3$d nonlinear elasticity, also in presence of
prestrain which is effectively manifested through the constraints of the
type (\ref{due.2}), can be found in the review paper \cite{LP} and
references therein.
\end{remark}

\v
{\bf Acknowledgments.}
The first author was partially supported by NSF  grant DMS-1714237, ``Models of controlled biological growth".
The second author was partially supported by NSF grants  DMS-1406730 and DMS-1613153. 

\v

\end{document}